\documentclass[11pt]{article}
\usepackage{amssymb,amsbsy,amsmath,amsfonts,amssymb,amscd,colordvi}
\usepackage{epsfig}
\usepackage{graphicx}
\usepackage[english]{babel}
\usepackage{times}
\usepackage{euscript}

\usepackage{tikz}
\usepackage{float}

\def \ve{\varepsilon}

\def \liminf{\mathop{\underline{\rm lim}}}
\def \be{\begin{equation}}
\def \ee{\end{equation}}

\def \eqdef{\mathop{=}\limits^{\hbox{\tiny def}}}

\usepackage{color}

\def\mx{\mathsf m}
\def\fr{\mathsf f}
\newtheorem{theorem}{Theorem}[section]
\newtheorem{lemma}{Lemma}[section]
\newtheorem{definition}[theorem]{Definition}
\newtheorem{remark}{\rm Remark\/}
\newtheorem{corollary}[theorem]{Corollary}
\newtheorem{proposition}[theorem]{Proposition}
\def \trait (#1) (#2) (#3){\vrule width #1pt height #2pt depth #3pt}
\def \fin{\hfill
	\trait (0.1) (5) (0)
	\trait (5) (0.1) (0)
	\kern-5pt
	\trait (5) (5) (-4.9)
	\trait (0.1) (5) (0)
\medskip}
\font \twbbb= msbm10 scaled \magstep0                 
\font \tenbbb= msbm7 scaled \magstep0                 
\newfam\bbbfam             
\textfont\bbbfam=\twbbb \scriptfont\bbbfam=\tenbbb    %

\begin{document}

\renewcommand{\theequation}{\thesection .\arabic{equation}}

\renewcommand{\baselinestretch}{1.2}

\title{Some remarks on the homogenization of immiscible incompressible two-phase flow in double porosity media}
\author{B. Amaziane$^{1}$, M. Jurak$^{2,}\footnote{Corresponding author.}$\ , L. Pankratov$^{1,3}$, A. Vrba\v ski$^4$}
\maketitle

\begin{small}
\noindent $^1$ Laboratoire de Math\'ematiques et de leurs Applications, CNRS-UMR 5142 Universit\'e de Pau,
Av. de l'Universit\'e, 64000 Pau, France.
E-mail: {\tt brahim.amaziane@univ-pau.fr}

\noindent $^2$ Faculty of Science, University of Zagreb, Bijeni\v cka 30, 10000 Zagreb, Croatia.
E-mail: {\tt jurak@math.hr}

\noindent $^3$ Laboratory of Fluid Dynamics and Seismic ({\sl RAEP 5top100}), Moscow Institute
of Physics and Technology, 9 Institutskiy per., Dolgoprudny, Moscow Region, 141700, Russian Federation.
E-mail: {\tt leonid.pankratov@univ-pau.fr}

\noindent $^4$  Faculty of Mining, Geology and Petroleum Engineering, University of Zagreb,
Pierottijeva 6, 10000 Zagreb, Croatia.
E-mail: {\tt anja.vrbaski@rgn.hr}
\end{small}

\begin{abstract}
This paper presents a study of immiscible incompressible two-phase flow through fractured porous media.
The results obtained earlier in the pioneer work by A. Bourgeat, S. Luckhaus,
A. Mikeli\'c (1996) and L.~M.~Yeh (2006) are revisited.
The main goal is to incorporate some of the most recent improvements in the convergence of the solutions in the homogenization of such models.
The microscopic model consists of the usual equations derived from the mass conservation of both fluids along with the Darcy-Muskat law.
The problem is written in terms of the phase formulation,
i.e. the saturation of one phase and the pressure of the second phase are primary unknowns.
We will consider a domain made up of several zones with different characteristics: porosity, absolute permeability,
relative permeabilities and capillary pressure curves.
The fractured medium consists of periodically repeating homogeneous blocks and
fractures, the permeability being highly discontinuous.
Over the matrix domain, the permeability is scaled by $\ve^\theta$, where $\ve$ is the size of a typical porous block and $\theta>0$ is a parameter.
The model involves highly oscillatory characteristics and internal nonlinear interface conditions.
Under some realistic assumptions on the data, the convergence of the solutions, and the macroscopic models corresponding to various
range of contrast are constructed using the two-scale convergence method combined with the dilation technique.
The results improve upon previously derived effective models to highly heterogeneous porous media with discontinuous capillary pressures.
\end{abstract}



\section{Introduction}
\label{intrrr}

\setcounter{equation}{0}

The modeling of displacement process involving two immiscible fluids in fractured porous media  is
important to many practical problems, including those in petroleum reservoir engineering, unsaturated zone hydrology,
and soil science. More recently, modeling multiphase flow has received an increasing attention in connection
with the disposal of radioactive waste and sequestration of $CO_2$. Furthermore, fractured rock domains
corresponding to the so-called Excavation Damaged Zone (EDZ) receives increasing attention in connection
with the behaviour of geological isolation of radioactive waste after the drilling of the wells or shafts, see, e.g., \cite{Shaw}.

In this paper we use the homogenization theory to derive a double porosity model
describing the flow of incompressible fluids in fractured reservoirs.
The model corresponds physically to immiscible incompressible two-phase flow through
fractured porous media. Naturally fractured reservoirs can be modeled by two
superimposed continua, a connected fracture system and a system of topologically
disconnected matrix blocks. The fracture system has  low storage capacity but high
conductivity, while the matrix block system has low conductivity and large storage capacity.
The majority of fluid transport will occur along flow paths through the fissure system.
When the system of fissures is so well developed that the matrix is
broken into individual blocks or cells that are isolated from each other, there is consequently
no flow directly from cell to cell, but only an exchange of fluid between each cell and the
surrounding fissure system. For more details on the physical formulation of such
problems see, e.g., \cite{bear,panf,van}.

This paper continues the research published in \cite{blm} and \cite{yeh2}, and the goal is to
reformulate in a more systematic manner and in somewhat more general context the homogenization
problem for an immiscible incompressible two-phase flow in double porosity media by weakening
the standing assumptions. Special attention is paid to developing a general approach to incorporating
highly heterogeneous porous media with discontinuous capillary pressures.

During recent decades mathematical analysis and homogenization of multiphase flows in porous media
have been the subject of investigation of many researchers owing to important applications in reservoir simulation.
There is an extensive literature on this subject. We will not attempt a literature review here but will
merely mention a few references. A recent review of the mathematical homogenization methods developed
for incompressible immiscible two-phase flow in porous media and compressible miscible flow in porous
media can be viewed in \cite{our-siam,HOS-2012,hor}.

Let us now turn to a brief review of the homogenization in double porosity media.
Here we restrict ourself to the mathematical homogenization method as described in \cite{hor} for flow and transport in porous media.
The interest for double porosity systems came at first from geophysics. The notion of double
porosity, or double permeability is borne from studies carried out on naturally fractured porous
rocks, such as oil fields. The double porosity model was
first introduced in \cite{bar} and it is since used in a wide range of
engineering specialities. The first rigorous mathematical result on the subject was obtained
in \cite{adh}, where a linear parabolic equation
with asymptotically degenerating coefficients describing a single-phase flow in
fractured media was considered.
This result is then generalized in \cite{bmp,bgpp,mk,sandr} for non-periodic
domains and various rates of contrast.
Linear double porosity models with thin fissures were considered in \cite{pr-AA,pr}.
A singular double porosity model was considered in \cite{bcp}. Notice that the works \cite{pr-AA,bgpp,mk,pr}
are done in the framework of Khruslov's energy characteristic method which is close to the $\Gamma$-convergence method. Let us also notice that the
double porosity model was obtained in \cite{hor} (see Chapter 3) using the two-scale convergence method.
Non-linear double porosity models, elliptic and parabolic,
including the homogenization in variable Sobolev spaces, were obtained in
\cite{edinburg-app,na-rwa-app,choq,edinb-catho,clark}. A study of discrete
double-porosity models in the case of elastic energies has been recently done in
\cite{braides-valer-piat}. Finally, in order to complete this brief review, we turn
to the multiphase flow double porosity models. These models were obtained e.g., in
\cite{blm,choq,yeh2} (see also \cite{hor} and the references therein) and recently in
\cite{latifa,ba-lp-doubpor} for immiscible compressible two-phase flows. A fully homogenized
model for incompressible two-phase flow in double porosity media was obtained in
\cite{Jurak}.

This paper is concerned with a nonlinear degenerate system of diffusion-convection equations  modeling the
flow and transport of immiscible incompressible fluids through highly heterogeneous porous media,
capillary and gravity effects being taken into account. We will consider a domain made up of several zones
with different characteristics: porosity, absolute permeability, relative permeabilities and capillary pressure curves.
The model to be presented herein is formulated in terms of the wetting phase saturation   and the non-wetting phase pressure,
and the feature of the global pressure as introduced in \cite{ant-kaz-mon1990,GC-JJ}
for incompressible immiscible flows is used to establish a priori estimates.
The governing equations are derived from the mass conservation laws of both fluids, along with constitutive
relations relating the velocities to the pressures gradients and gravitational effects. Traditionally, the standard Muskat-Darcy law provides this relationship.
Let us mention that the main difficulties related to the mathematical analysis of such equations
are the coupling and the degeneracy of the diffusion term in the saturation equation.
Moreover the transmission conditions are nonlinear and the saturation is discontinuous at the interface separating the two media.

We start with a microscopic model defined on a domain with periodic microstructure.
We will consider a domain made up of several zones with different characteristics: porosity, absolute permeability, relative permeabilities and capillary pressure curves.
The fractured medium consists of periodically repeating homogeneous blocks and
fractures, the permeability being highly discontinuous.
Over the matrix domain, the permeability is scaled by $\ve^\theta$, where $\ve$ is the size of a typical porous block and $\theta>0$ is a parameter.
Our aim is to study the macroscopic behavior of solutions of this system of equations as $\ve$ tends
to zero and give a rigorous mathematical derivation of upscaled models by means of the two-scale convergence method combined with the dilation technique.
Thus, we extend the results of \cite{blm,yeh2} to the case of highly heterogeneous porous media with discontinuous capillary pressures.

The rest of the paper is organized as follows. In Section \ref{micromodel},
we describe the physical model and formulate the corresponding mathematical problem.
We also provide the assumptions on the data and a weak formulation of the problem firstly in terms
of phase pressures and secondly in terms of the global pressure and the saturation.
Section \ref{uni-est} is then devoted to the derivation of the basic {\it a priori} estimates
of the problem under consideration. In Section \ref{ex-comp-sf} we formulate the two-scale
convergence results which will be used in the derivation of the homogenized system.
The key point here is the proof of the compactness result for the restriction-extension
sequence of the wetting fluid saturation defined on the
fracture set. It is done by using the ideas from \cite{yeh2}. Section \ref{dil-oper} is devoted
to the definition and the properties of the dilation operator and to the formulation
of the convergence results for the dilated functions defined on the matrix part. The key point
of the section is the proof of the compactness result for the dilated saturations
which is done by using the compactness result from
\cite{our-siam}. The formulation of the main results of the paper is given in
Section \ref{main-res}.
The resulting homogenized problem is a dual-porosity type model that contains a term representing
memory effects which could be seen as source term or as a time delay for $\theta = 2$, and it is a single
porosity model with effective coefficients for $0 < \theta < 2$ or $\theta > 2$.
The proof of the convergence theorem in the critical case ($\theta =2$) is
done in subsection \ref{hig-crtic-subsec}. The key point here is subsection \ref{identif-l-theta}, where
we prove the uniqueness of the solution to the local problem. The proof is done by reducing
the problem in the phase formulation to a boundary value problem for an imbibition equation
and by using ideas from \cite{vazquez}. The proofs of the convergence theorems for non-critical
cases ($\theta>2$ or $0<\theta<2$) are given in subsections \ref{very-hig-crtic-subsec}, \ref{moder-case-subsec}.
The effective model obtained in the case of moderate contrast ($0<\theta<2$, subsection \ref{moder-case-subsec}),
up to our knowledge, is for the first time proposed and rigorously justified here.

\section{Formulation of the  problem}
\label{micromodel}
\setcounter{equation}{0}

The outline of this section is as follows. First, in subsection \ref{micromodel-simpl}
we give a short description of the mathematical and physical model used in this study
for immiscible incompressible two-phase flow in a periodic double porosity medium.
The notion of the global pressure is briefly recalled in subsection \ref{gp-relat}.
Finally, in subsection \ref{def-weak-sol}, we present the main assumptions on the data and we define the weak solution to our problem, first in terms of phase pressures and then an equivalent one in terms of the global pressure and saturation.

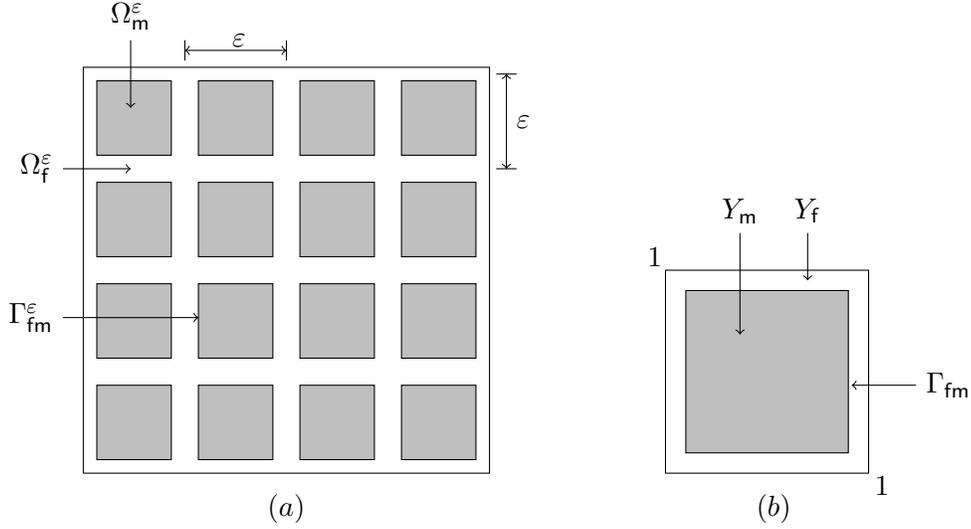
\begin{figure}[H]
		\begin{center}
			\begin{tikzpicture}[scale=0.9]
				\def\rectanglepath{-- ++(1.1cm,0cm)  -- ++(0cm,1.1cm) -- ++(-1.1cm,0cm) -- cycle}
				
				\filldraw[fill=white!30] (-0.2,-0.2) -- (5.8,-0.2) --  (5.8, 5.8) -- (-0.2, 5.8) -- cycle;
				
				\filldraw[fill=lightgray] (0,  0) \rectanglepath;
				\filldraw[fill=lightgray] (1.5,0) \rectanglepath;
				\filldraw[fill=lightgray] (3,  0) \rectanglepath;
				\filldraw[fill=lightgray] (4.5,0) \rectanglepath;
				
				\filldraw[fill=lightgray] (0,  1.5) \rectanglepath;
				\filldraw[fill=lightgray] (1.5,1.5) \rectanglepath;
				\filldraw[fill=lightgray] (3,  1.5) \rectanglepath;
				\filldraw[fill=lightgray] (4.5,1.5) \rectanglepath;
				
				\filldraw[fill=lightgray] (0,  3.0) \rectanglepath;
				\filldraw[fill=lightgray] (1.5,3.0) \rectanglepath;
				\filldraw[fill=lightgray] (3,  3.0) \rectanglepath;
				\filldraw[fill=lightgray] (4.5,3.0) \rectanglepath;
				
				\filldraw[fill=lightgray] (0,  4.5) \rectanglepath;
				\filldraw[fill=lightgray] (1.5,4.5) \rectanglepath;
				\filldraw[fill=lightgray] (3,  4.5) \rectanglepath;
				\filldraw[fill=lightgray] (4.5,4.5) \rectanglepath;
				
				\draw[<-] (0.5,5.2) -- (0.5,6.2) node[anchor=south] {$\Omega_\mx^{\varepsilon}$};
				\draw[<-] (0.5,4.3) -- (-0.5,4.3) node[anchor=east] {$\Omega_\fr^{\varepsilon}$};
				\draw[<-] (1.5,2.1) -- (-0.5,2.1) node[anchor=east] {$\Gamma^\ve_{\fr\mx}$};
				
				\draw (1.3,5.9)--(1.3,6.2);
				\draw (2.8,5.9)--(2.8,6.2);
				\draw[<->] (1.3,6.05)--(2.8,6.05);
				\draw (2.1,6.2) node {$\varepsilon$};
				
				\draw (5.9,5.7)--(6.2,5.7);
				\draw (5.9,4.3)--(6.2,4.3);
				\draw[<->] (6.05,4.3)--(6.05,5.7);
				\draw (6.3,5.0) node {$\varepsilon$};
						
				\node [below=0.5cm] at (2.8,0.2) {$(a)$};
						
				
				\def\rectanglepath{-- ++(1.1cm,0cm)  -- ++(0cm,1.1cm) -- ++(-1.1cm,0cm) -- cycle}
				
				\draw[fill=white!30] (8.4,-0.2) -- (11.4,-0.2) -- (11.4,2.8) -- (8.4,2.8) -- cycle;
				\draw[fill=lightgray]   (8.7,0.1) -- (11.1,0.1) -- (11.1,2.5) -- (8.7,2.5) -- cycle;
				
				\draw (11.6,-0.1) node[anchor=north] {$1$};
				\draw (8.5,3.0) node[anchor=east] {$1$};
				
				\draw[<-] (9.5,1.85) -- (9.5,3.35) node[anchor=south] {$Y_\mx$};
				\draw[<-] (10.5,2.65) -- (10.5,3.35) node[anchor=south] {$Y_\fr$};
				\draw[<-] (11.15,1.1) -- (12.1,1.1) node[anchor=west] {$\Gamma_{\fr\mx}$};
				
				\node [below=0.5cm] at (10.0,0.2) {$(b)$};
			\end{tikzpicture}
		\end{center}
		\caption{(a) The domain $\Omega$. \quad (b) The reference cell $Y$.}
		\label{fig:ref}
	\end{figure}

\subsection{Microscopic model}
\label{micromodel-simpl}

We consider a reservoir $\Omega \subset \mathbb{R}^d$ ($d = 2, 3$)
which is assumed to be a bounded, connected Lipschitz domain with a periodic microstructure.
More precisely, we will scale this structure by a
parameter $\ve$ which represents the ratio of the cell size to
the size of the whole region $\Omega$  and we assume that $0 < \ve \ll 1$ is a small parameter
tending to zero. Let $Y = (0, 1)^d$ be a basic cell of a fractured porous medium.
For the sake of simplicity and without loss of generality, we assume that $Y$ is made
up of two homogeneous porous media $Y_\mx$ and $Y_\fr$ corresponding to the parts
of the domain occupied by the matrix block and the fracture, respectively
(see Fig.\ref{fig:ref} (b)). Thus $Y = Y_\mx \cup Y_\fr \cup \Gamma_{\fr\mx}$,
where $\Gamma_{\fr\mx}$ denotes
the interface between the two media. Let $\Omega^\ve_\ell$ with $\ell = "\fr"$ or $"\mx"$
denotes the open set corresponding to the porous medium with index $\ell$. Then
$\Omega = \Omega^\ve_\mx \cup \Gamma^\ve_{\fr\mx} \cup \Omega^\ve_\fr$,
where $\Gamma^\ve_{\fr\mx} \eqdef \partial \Omega^\ve_\fr \cap \partial \Omega^\ve_\mx \cap \Omega$
and the subscripts $"\mx"$, $"\fr"$ refer to the matrix and fracture, respectively
(see Fig.\ref{fig:ref} (a)).
For the sake of simplicity, we assume that $\Omega^\ve_\mx \cap \partial \Omega =
\emptyset$. We also introduce the notation:
\begin{equation}
\label{omeg12++}
\Omega_T \eqdef \Omega \times (0,T), \quad
\Omega^\ve_{\ell,T} \eqdef \Omega^\ve_\ell \times (0,T), \quad
\Sigma^\ve_T \eqdef \Gamma^\ve_{\fr\mx} \times (0,T), \quad {\rm where}\,\, T > 0 \,\, {\rm is \,\, fixed.}
\end{equation}

Before describing the equations of the model, we give some notation:
$\Phi^\ve(x) = \Phi(x, \frac{x}{\ve})$ is the porosity of the reservoir $\Omega$;
$K^\ve(x) = K(x, \frac{x}{\ve})$ is the absolute permeability tensor of $\Omega$;
$\varrho_w$, $\varrho_n$ are the densities of wetting and nonwetting fluids,
respectively; $S^\ve_{\ell, w} = S^\ve_{\ell, w}(x, t)$, $S^\ve_{\ell, n}
= S^\ve_{\ell, n}(x, t)$ are the saturations of wetting and nonwetting fluids
in $\Omega^\ve_\ell$, respectively;
$k_{r,w}^{(\ell)} = k_{r,w}^{(\ell)}(S^\ve_{\ell, w})$,
$k_{r,n}^{(\ell)} = k_{r,n}^{(\ell)}(S^\ve_{\ell, n})$ are the
relative permeabilities of wetting and nonwetting fluids in
$\Omega^\ve_\ell$, respectively; $p^\ve_{\ell, w} =
p^\ve_{\ell, w}(x,t)$, $p^\ve_{\ell, n} = p^\ve_{\ell, n}(x,t)$ are the
phase pressures of wetting and nonwetting fluids in $\Omega^\ve_\ell$,
respectively. Here $\ell = \fr, \mx$.

The conservation of mass in each phase can be written as
(see, e.g., \cite{GC-JJ,ZC-GH-YM-06,HR}):
\begin{equation}
\label{debut1}
\left\{
\begin{array}[c]{ll}
\displaystyle
\Phi^\ve(x) \frac{\partial}{\partial t}\left[S^\ve_{\ell, w}\,\varrho_w(p^\ve_{\ell, w})\right] +
{\rm div} \big\{\varrho_w(p^\ve_{\ell, w}) \, \vec q^{\,\ve}_{\ell, w} \big\} = F^\ve_{\ell,w}(x,t)
\quad {\rm in}\,\, \Omega^\ve_{\ell,T}; \\[3mm]
\displaystyle
\Phi(x) \frac{\partial}{\partial t} \left[S^\ve_{\ell, n}\, \varrho_n(p^\ve_{\ell, n})\right] +
{\rm div}\big\{\varrho_n(p^\ve_{\ell, n}) \, \vec q^{\,\ve}_{\ell, n} \big\} = F^\ve_{\ell,n}(x,t)
\quad {\rm in}\,\, \Omega^\ve_{\ell,T},  \\[3mm]
\end{array}
\right.
\end{equation}
where the velocities of the wetting and nonwetting fluids $\vec q^{\,\ve}_{\ell, w}$,
$\vec q^{\,\ve}_{\ell, n}$ are defined by Darcy-Muskat's law:
\begin{equation}
\vec q^{\,\ve}_{\ell, w} \eqdef -K^\ve(x) \lambda_{\ell, w}(S^\ve_{\ell, w})
\left[\nabla p^\ve_{\ell, w} - \varrho_w(p^\ve_{\ell, w})\, \vec{g}\right],
\quad \!{\rm with}\,\, \lambda_{\ell, w}(S^\ve_{\ell, w}) \eqdef
\frac{k_{r,w}^{(\ell)}}{\mu_{w}}(S^\ve_{\ell, w});
\label{eq.qw}
\end{equation}
\begin{equation}
\vec q^{\,\ve}_{\ell, n} \eqdef - K^\ve(x) \lambda_{\ell, n}(S^\ve_{\ell, n})
\left[\nabla p^\ve_{\ell, n} - \varrho_n(p^\ve_{\ell, n})\, \vec{g}\right], \quad
{\rm with}\,\, \lambda_{\ell, n}(S_{\ell, n}) \eqdef \frac{k_{r,n}^{(\ell)}}{\mu_{n}}
(S^\ve_{\ell, n}).
\label{eq.qg}
\end{equation}
Here $\vec g$, $\mu_w, \mu_n$ are the gravity vector and the viscosities of the
wetting and nonwetting fluids, respectively. The source terms $F^\ve_{\ell,w}, F^\ve_{\ell,n}$ are
given by:
\begin{equation}
\label{sour1}
F^\ve_{\ell,w} \eqdef \varrho_w(p^\ve_{\ell, w}) S^I_{\ell, w} f_I(x,t) -
\varrho_w(p^\ve_{\ell, w}) S^\ve_{\ell, w} f_P(x,t);
\end{equation}
\begin{equation}
\label{sour2}
F^\ve_{\ell,n} \eqdef \varrho_n(p^\ve_{\ell, n}) S^I_{\ell, n} f_I(x,t) -
\varrho_n(p^\ve_{\ell, n}) S^\ve_{\ell, n} f_P(x,t),
\end{equation}
where $f_I, f_P \geqslant 0$ are injection and productions terms and
$S^I_{\ell, w}, S^I_{\ell, n}$ are known injection saturations.

From now on we deal with two incompressible fluids, that is the densities of the wetting
and nonwetting fluids are constants, which for the sake of simplicity and brevity, will
be taken equal to one, i.e. $\varrho_w(p^\ve_{\ell,w}) = \varrho_n(p^\ve_{\ell,n}) = 1$.
The model is completed as follows. By the definition of saturations, one has
$S^\ve_{\ell, w} + S^\ve_{\ell, n} = 1$ with $S^\ve_{\ell, w}, S^\ve_{\ell, n} \geqslant 0$.
We set $S^\ve_\ell \eqdef S^\ve_{\ell, w}$.
Then the curvature of the contact surface between the two fluids links the difference in
the pressures of the two phases to the saturation by the capillary pressure law:
\begin{equation}
P_{\ell,c}(S^\ve_\ell) \eqdef p^\ve_{\ell, n} - p^\ve_{\ell, w} \quad {\rm with} \,\,
P^\prime_{\ell,c}(s) < 0\,\, {\rm for\,\, all}\,\,s \in (0, 1)
\,\, {\rm and} \,\, P_{\ell,c}(1) = 0,
\label{eq.pc}
\end{equation}
where $P^\prime_{\ell,c}(s)$ denotes the derivative of the function $P_{\ell,c}(s)$.

Now due to the assumptions on the densities of the liquids, we rewrite the system
(\ref{debut1}) as follows:
\begin{equation}
\label{debut2}
\left\{
\begin{array}[c]{ll}
0 \leqslant {S}^\ve \leqslant 1 \quad {\rm in}\,\, \Omega_T;\\[3mm]
\displaystyle
\Phi^\ve(x) \frac{\partial {S}^\ve}{\partial t} - {\rm div}\, \left\{K^\ve(x)
\lambda_{w}\left(\frac{x}{\ve}, {S}^\ve\right)
\left(\nabla {\mathsf p}^\ve_{w} - \vec{g}\right)\right\} = F^\ve_{w} \quad {\rm in}\,\, \Omega_{T}; \\[4mm]
\displaystyle
-\Phi^\ve(x) \frac{\partial {S}^\ve}{\partial t} -
{\rm div}\, \left\{K^\ve(x) \lambda_{n} \left(\frac{x}{\ve}, {S}^\ve\right)
\left(\nabla {\mathsf p}^\ve_{n} - \vec{g} \right)\right\} =  F^\ve_{n} \quad {\rm in}\,\, \Omega_{T}; \\[4mm]
P_{c}\left(\frac{x}{\ve}, {S}^\ve\right) = {\mathsf p}^\ve_{n} - {\mathsf p}^\ve_{w}
\quad {\rm in}\,\, \Omega_{T},
\end{array}
\right.
\end{equation}
where $\lambda_{\ell,n}(S^\ve_\ell) := \lambda_{\ell,n}(1-S^\ve_\ell)$ and each function
$u^\ve := {S}^\ve, {\mathsf p}^\ve_{w}, {\mathsf p}^\ve_{n},  F^\ve_{w},  F^\ve_{n}$
is defined as:
\begin{equation}
\label{nota-pwat-pgaz}
u^\ve \eqdef  u^\ve_{\fr}(x, t)\, {\bf 1}^\ve_\fr(x) +  u^\ve_{\mx}(x, t)\, {\bf 1}^\ve_\mx(x).
\end{equation}
Here ${\bf 1}^\ve_\ell = {\bf 1}_\ell(\frac{x}{\ve})$ is the characteristic function of the subdomain
$\Omega^\ve_\ell$ for $\ell = \fr, \mx$. The exact form of the porosity function and the absolute
permeability tensor corresponding to the double porosity model will be specified in conditions
{(A.1)}, {(A.2)} in subsection \ref{def-weak-sol} below.

Model (\ref{debut2}) have to be completed with appropriate interface, boundary
and initial conditions.

\noindent Interface conditions.
The continuity at the interface $\Gamma^\ve_{\fr\mx}$ of the
phase fluxes and the phase pressures, gives the following transmission conditions:
\begin{equation}
\label{inter-condit}
\left\{
\begin{array}[c]{ll}
\vec  q^{\,\ve}_{\fr,w} \cdot \vec \nu = \vec  q^{\,\ve}_{\mx,w}\cdot \vec \nu \quad {\rm and} \quad
\vec  q^{\,\ve}_{\fr,n} \cdot \vec \nu = \vec  q^{\,\ve}_{\mx,n}\cdot \vec \nu
\quad {\rm on}\,\, \Sigma^\ve_T; \\[2mm]
p^\ve_{\fr,w} = p^\ve_{\mx,w} \quad {\rm and} \quad p^\ve_{\fr,n} = p^\ve_{\mx,n} \quad {\rm on}\,\, \Sigma^\ve_T, \\
\end{array}
\right.
\end{equation}
where $\Sigma^\ve_T$ is defined in (\ref{omeg12++}),
$\vec \nu$ is the unit outer normal on $\Gamma^\ve_{\fr\mx}$, and the fluxes
$\vec q^{\,\ve}_{\ell,w}, \vec q^{\,\ve}_{\ell,n}$, under the assumption on the densities of the liquids,
are equal to the velocities (\ref{eq.qw}), (\ref{eq.qg}).

\begin{remark}
\label{interface-saturrr}
It is important to notice that in contrast to the functions
${\mathsf p}^\ve_n, {\mathsf p}^\ve_w$, the saturation ${S}^\ve$
may have a jump at the interface $\Gamma^\ve_{\fr\mx}$. Namely, it is easy to
see from the transmission conditions (\ref{inter-condit}) for the phase pressures
that $P_{\fr,c}(S^\ve_1) = P_{\mx,c}(S^\ve_2)$ on $\Sigma^\ve_T$ which gives a discontinuity of
the saturation at the interface.
\end{remark}

Now we specify the boundary and initial conditions. We suppose that the boundary
$\partial \Omega$ consists of two parts $\Gamma_{1}$ and $\Gamma_{2}$ such
that $\Gamma_{1} \cap \Gamma_{2} = \emptyset$, $\partial \Omega =
\overline\Gamma_{1} \cup \overline\Gamma_{2}$.

\noindent Boundary conditions:
\begin{equation}
\label{bc3}
\left\{
\begin{array}[c]{ll}
{\mathsf p}^\ve_{w}(x, t) = {\mathsf p}^\ve_{n}(x, t) = 0 \quad {\rm on} \,\,
\Gamma_{1} \times (0,T); \\[2mm]
\vec q^{\,\ve}_{\fr,w} \cdot \vec \nu = \vec q^{\,\ve}_{\fr,n} \cdot \vec \nu = 0
\quad {\rm on} \,\, \Gamma_{2} \times (0,T).\\
\end{array}
\right.
\end{equation}
\noindent Initial conditions:
\begin{equation}
\label{init1}
{\mathsf p}^\ve_{w}(x, 0) = {\mathsf p}_{w}^{\bf 0}(x) \quad {\rm and}
\quad
{\mathsf p}^\ve_{n}(x, 0) = {\mathsf p}_{n}^{\bf 0}(x)
\quad {\rm in} \,\, \Omega.
\end{equation}

\subsection{A fractional flow formulation}
\label{gp-relat}

In the sequel, we will use a formulation obtained after transformation using
the concept of the global pressure introduced in \cite{ant-kaz-mon1990,GC-JJ}.
For each subdomain $\Omega^\ve_\ell$, the global pressure, ${\mathsf P}^\ve_\ell$,
is defined by:
\begin{equation}
\label{gp1}
p^\ve_{\ell,w} \eqdef {\mathsf P}^\ve_\ell + {\mathsf G}_{\ell,w}(S^\ve_\ell)
\quad {\rm and} \quad
p^\ve_{\ell,n} \eqdef {\mathsf P}^\ve_\ell + {\mathsf G}_{\ell,n}(S^\ve_\ell),
\end{equation}
where the functions ${\mathsf G}_{\ell,w}(s)$, ${\mathsf G}_{\ell,n}(s)$
are given by:
\begin{equation}
\label{gp3}
{\mathsf G}_{\ell,n}(S^\ve_\ell) \eqdef
{\mathsf G}_{\ell,n}(0) + \int\limits_0^{S^\ve_\ell} \frac{\lambda_{\ell,w}(s)}
{\lambda_{\ell}(s)} \,P^\prime_{\ell,c}(s)\, ds\quad {\rm and} \quad
{\mathsf G}_{\ell, w}(S^\ve_\ell) \eqdef {\mathsf G}_{\ell,n}(S^\ve_\ell)
- P_{\ell,c}\left(S^\ve_\ell\right),
\end{equation}
where $\lambda_{\ell}(s) \eqdef \lambda_{\ell,w}(s) + \lambda_{\ell,n}(s)$ and ${\mathsf G}_{\ell,n}(0)$ is a constant
chosen to ensure $p^\ve_{\ell,w}  \leqslant {\mathsf P}^\ve_\ell \leqslant p^\ve_{\ell,n}$.
Notice that from (\ref{gp3}) we get:
\begin{equation}
\label{gp+grad}
\lambda_{\ell,w}(S^\ve_\ell) \nabla {\mathsf G}_{\ell,w}(S^\ve_\ell) = \nabla\beta_\ell(S^\ve_\ell)
\quad {\rm and} \quad
\lambda_{\ell,n}(S^\ve_\ell) \nabla {\mathsf G}_{\ell,n}(S^\ve_\ell) = - \nabla\beta_\ell(S^\ve_\ell),
\end{equation}
where
\begin{equation}
\label{upsi-1}
\beta_\ell(s) \eqdef \int\limits_0^s \alpha_\ell(\xi)\, d\xi \quad {\rm with} \,\,\,
\alpha_\ell(s) \eqdef
\frac{\lambda_{\ell,n}(s)\, \lambda_{\ell,w}(s)} {\lambda_\ell(s)} \left| P^\prime_{\ell,c}(s) \right|.
\end{equation}

Furthermore, we have the following important relation:
\begin{equation}
\label{gp6-new}
\lambda_{\ell,n} (S^\ve_\ell) |\nabla p^\ve_{\ell,n}|^2 +
\lambda_{\ell,w} (S^\ve_\ell) |\nabla p^\ve_{\ell,w}|^2
=
\lambda_{\ell} (S^\ve_\ell) |\nabla {\mathsf P}^\ve_\ell |^2 +
\left|\nabla \mathfrak{b}_\ell(S^\ve_\ell) \right|^2,
\end{equation}
where
\begin{equation}
\label{bbb}
\mathfrak{b}_\ell(s) \eqdef \int\limits_0^s \mathfrak{a}_\ell(\xi)\, d\xi
\quad {\rm with} \,\, \mathfrak{a}_\ell(s) \eqdef
\sqrt{\frac{\lambda_{\ell,n}(s)\, \lambda_{\ell,w}(s)} {\lambda_\ell(s)}}\, \left| P^\prime_{\ell,c}(s) \right|.
\end{equation}

Now if we use the global pressure and the saturation as new unknown
functions then (\ref{debut2}) reads:
\begin{equation}
\label{gp+3-simp}
\left\{
\begin{array}[c]{ll}
0 \leqslant S^\ve_\ell \leqslant 1 \quad {\rm in}\,\, \Omega^\ve_{\ell,T};\\[3mm]
\displaystyle
\Phi^\ve(x) \frac{\partial S^\ve_\ell}{\partial t} -
{\rm div}\, \bigg\{K^\ve(x) \left[
\lambda_{\ell,w} (S^\ve_\ell) \nabla {\mathsf P}^\ve_\ell +
\nabla \beta_\ell(S^\ve_\ell) - \lambda_{\ell,w} (S^\ve_\ell) \vec{g} \right] \bigg\}
=  F^\ve_{\ell,w} \quad {\rm in}\,\, \Omega^\ve_{\ell,T}; \\[4mm]
\displaystyle
-\Phi^\ve(x) \frac{\partial S^\ve_\ell}{\partial t}
-
{\rm div}\, \bigg\{K^\ve(x) \left[ \lambda_{\ell,n} (S^\ve_\ell)
\nabla {\mathsf P}^\ve_\ell - \nabla \beta_\ell(S^\ve_\ell)
-
\lambda_{\ell,n}(S^\ve_\ell) \vec{g} \right] \bigg\} =
F^\ve_{\ell,n}\quad {\rm in}\,\, \Omega^\ve_{\ell,T}.\\
\end{array}
\right.
\end{equation}

The system (\ref{gp+3-simp}) is completed by the following boundary, interface and initial conditions.

\noindent Boundary conditions:
\begin{equation}
\label{bc3+gp}
\left\{
\begin{array}[c]{ll}
{S}^\ve = 1 \,\,\,{\rm and}\,\,\,
{\mathsf P}^\ve = {\mathsf P}_{\Gamma_1} \quad {\rm on} \,\,
\Gamma_{1} \times (0,T); \\[2mm]
\vec q^{\,\ve}_{\fr,w} \cdot \vec \nu = \vec q^{\,\ve}_{\fr,n} \cdot \vec \nu = 0
\quad {\rm on} \,\, \Gamma_{2} \times (0,T),\\
\end{array}
\right.
\end{equation}
where ${\mathsf P}_{\Gamma_1}$ is a given constant and $\vec q^{\,\ve}_{\ell, w}, \vec q^{\,\ve}_{\ell, n}$ are defined by
\begin{equation}
\label{fluxes+gp-w}
\vec q^{\,\ve}_{\ell, w} \eqdef -K^\ve(x) \left[
\lambda_{\ell,w} (S^\ve_\ell) \nabla {\mathsf P}^\ve_\ell +
\nabla \beta_\ell(S^\ve_\ell) - \lambda_{\ell,w} (S^\ve_\ell) \vec{g} \right];
\end{equation}
\begin{equation}
\label{fluxes+gp-n}
\vec q^{\,\ve}_{\ell, n} \eqdef - K^\ve(x) \left[ \lambda_{\ell,n} (S^\ve_\ell)
\nabla {\mathsf P}^\ve_\ell - \nabla \beta_\ell(S^\ve_\ell)
- \lambda_{\ell,n}(S^\ve_\ell) \vec{g} \right].
\end{equation}

\noindent Interface conditions:
\begin{equation}
\label{ic-glob}
\left\{
\begin{array}[c]{ll}
\vec q^{\,\,\ve}_{\fr,w} \cdot \vec \nu =
\vec q^{\,\,\ve}_{\mx,w}\cdot \vec \nu \,\,\, {\rm and} \,\,\,
\vec q^{\,\,\ve}_{\fr,n} \cdot \vec \nu =
\vec q^{\,\,\ve}_{\mx,n}\cdot \vec \nu
\quad {\rm on}\,\, \Sigma^\ve_T; \\[3mm]
{\mathsf P}^\ve_\fr + {\mathsf G}_{\fr,j}(S^\ve_\fr) =
{\mathsf P}^\ve_\mx + {\mathsf G}_{\mx,j}(S^\ve_\mx)
\quad {\rm on}\,\, \Sigma^\ve_T \quad (j = w, n);\\[3mm]
P_{\fr,c}\left(S^\ve_\fr\right) = P_{\mx,c}\left(S^\ve_\mx\right)
\quad {\rm on}\,\, \Sigma^\ve_T.
\end{array}
\right.
\end{equation}

Note that the global pressure function might be discontinuous
at the interface. This makes the compactness
result in Section \ref{ex-comp-sf} non-trivial.

\noindent Initial conditions:
\begin{equation}
\label{init1-gl+gp}
S^\ve_\ell(x, 0) = S^{\bf 0}_\ell(x) \,\, {\rm and} \,\,
{\mathsf P}^\ve_\ell(x, 0) = {\mathsf P}^{\bf 0}_\ell(x) \quad {\rm in} \,\, \Omega.
\end{equation}

\subsection{Weak formulations of the problem}
\label{def-weak-sol}

Let us begin this subsection by stating the following assumptions.

\begin{itemize}

\item[{(A.1)}] The porosity $\Phi^\ve$ is given by
$
\Phi^\ve(x) \eqdef \Phi^\ve_{\fr}(x)\, {\bf 1}^\ve_\fr(x) +  \Phi^\ve_{\mx}(x)\, {\bf 1}^\ve_\mx(x) =
\Phi^\ve_{\fr}(x)\, {\bf 1}^\ve_\fr(x) +  \Phi_{\mx}\left(\frac{x}{\ve}\right)\, {\bf 1}^\ve_\mx(x),
$
where $\Phi^\ve_{\fr} \in L^{\infty}(\Omega)$ and
there are positive constants $0 < \phi^{\ell}_- < \phi^{\ell}_+ < 1$, $\ell = \fr, \mx$, that do not depend on $\ve$ and such that
$0 <  \phi^{\fr}_- \leqslant \Phi^\ve_{\fr}(x) \leqslant \phi^{\fr}_+ < 1$ a.e. in $\Omega$.
Moreover, $\Phi^\ve_{\fr} \longrightarrow \Phi^{\rm H}_{\fr}$
strongly in $L^2(\Omega)$.
$\Phi_\mx = \Phi_\mx(y)$ is $Y$-periodic, $\Phi_{\mx} \in L^{\infty}(Y)$
and such that
$0 <  \phi^{\mx}_- \leqslant \Phi_\mx(y) \leqslant \phi^{\mx}_+ < 1$ a.e. in $Y$.

\item[{(A.2)}] The permeability $K^\ve(x) = K^\ve(x, \frac{x}{\ve})$ is defined as
$
\label{tensork-expr}
K^\ve(x, y) \eqdef K(x, y)\, {\bf 1}^\ve_\fr(x) + \varkappa(\ve)\,  K(x, y)\, {\bf 1}^\ve_\mx(x),
$
where $\varkappa(\ve) \eqdef \ve^\theta$ with $\theta > 0$ and $K \in (L^{\infty}(\Omega\times Y))^{d\times d}$.
Moreover, there exist constants $k_{\rm min}, k^{\rm max}$ such that  $0 < k_{\rm min} < k^{\rm max}$ and
$
\label{tensork}
k_{\rm min} |\xi|^2 \leq (K(x, y)\,\xi, \xi) \leq k^{\rm max} |\xi|^2 \,\, {\rm for \, all \,
\xi \in \mathbb{R}^d, \,\, a.e. \, in}\,\, \Omega\times Y.
$

\item[{(A.3)}] The capillary pressure function $P_{\ell,c}(s) \in C^1([0, 1]; \mathbb{R}^+)$, $\ell = \fr, \mx$.
Moreover, $P_{\ell,c}^\prime(s) < 0$ in $[0, 1]$, $P_{\ell,c}(1) = 0$  and $P_{\fr,c}(0) = P_{\mx,c}(0)$.

\item[{(A.4)}] The functions $\lambda_{\ell,w}, \lambda_{\ell,n}$
belong to the space $C([0, 1]; \mathbb{R}^+)$ and satisfy the following properties: \\
{(i)} $0 \leqslant \lambda_{\ell,w}, \lambda_{\ell,n} \leqslant 1$ in $[0, 1]$;
{(ii)} $\lambda_{\ell,w}(0) = 0$ and $\lambda_{\ell,n}(1) = 0$; {(iii)} there
is a positive constant $L_0$ such that $\lambda_{\ell}(s) =  \lambda_{\ell,w}(s) + \lambda_{\ell,n}(s)
\geqslant L_0 > 0$ in $[0, 1]$.

\item[{(A.5)}] The functions $\alpha_\ell \in C([0, 1]; \mathbb{R}^+)$. Moreover,
$\alpha_\ell (0) = \alpha_\ell (1) = 0$ and $\alpha_\ell > 0$ in $(0, 1)$.

\item[{(A.6)}] The functions $\beta^{-1}_\ell$, inverse of $\beta_\ell$ defined in (\ref{upsi-1})
are H\"older functions of order $\gamma\in (0, 1)$ in $[0, \beta_\ell(1)]$.
Namely, there exists a positive constant $C_\beta$ such that for all $s_1, s_2 \in [0, \beta(1)]$,
we have:
$$
\left|\beta^{-1}_\ell(s_1) - \beta^{-1}_\ell(s_2) \right| \leqslant C_\beta \, |s_1 - s_2|^\gamma.
$$

\item[{(A.7)}] The initial data for the pressures are such that ${\mathsf p}^{\bf 0}_{n},
{\mathsf p}^{\bf 0}_{w} \in L^2(\Omega)$.

\item[{(A.8)}] The initial data for the saturation $S^{\bf 0}$ is given by
$P_{\ell,c}(S^{\bf 0}_\ell) ={\mathsf p}^{\bf 0}_{\ell, n} - {\mathsf p}^{\bf 0}_{\ell, w}$
    and is such that
${S}^{\bf 0} \in L^\infty(\Omega)$ and $0 \leqslant {S}^{\bf 0} \leqslant 1$ a.e.in $\Omega$.

\item[{(A.9)}] The source terms $F^\ve_{w}, F^\ve_{n}$ are equal to zero on the matrix part, i.e.
$
\label{sour1-A9}
F^\ve_{w} \eqdef {\bf 1}^\ve_\fr(x)\, \big[S^I_{\fr, w} f_I(x,t) - S^\ve_{\fr} f_P(x,t)\big]
$
and
$
F^\ve_{n} \eqdef {\bf 1}^\ve_\fr(x)\, \big[S^I_{\fr, n} f_I(x,t) - (1 - S^\ve_{\fr}) f_P(x,t)\big],
$
where $f_I, f_P \in L^2(\Omega_T)$ and $0 \leqslant S^I_{\fr, w}, S^I_{\fr, n} \leqslant 1$.

\end{itemize}

The assumptions (A.1)--(A.9) are classical and physically meaningful for existence results and
homogenization problems of two-phase flow in porous media. They are similar to the assumptions made in \cite{ant-kaz-mon1990,GC-JJ}
that dealt with the existence of a weak solution of the studied problem.

We next introduce the following Sobolev space:
$
H^1_{\Gamma_{1}}(\Omega) \eqdef \left\{u \in H^1(\Omega) \,:\,
u = 0 \,\, {\rm on}\,\, \Gamma_{1} \right\}.
$
The space $H^1_{\Gamma_{1}}(\Omega)$ is a Hilbert space. The norm in this space
is given by $\Vert u \Vert_{H^1_{\Gamma_{1}}(\Omega)} = \Vert \nabla u \Vert_{(L^2(\Omega))^d}$.

\begin{definition}[\sl Weak solution in terms of phase pressures]
\label{def-weak-simp}

We say that the functions $\langle {\mathsf p}^\ve_{w}, {\mathsf p}^\ve_{n} , {S}^\ve\rangle$ is
a weak solution of {\rm problem (\ref{debut2})} if

\begin{itemize}

\item[{(i)}] $0 \leqslant {S}^\ve \leqslant 1$ a.e.\ in $\Omega_T$ and $P_{\ell,c}(S^\ve_\ell) \eqdef {\mathsf p}^\ve_{\ell, n} - {\mathsf p}^\ve_{\ell, w}$
    for $\ell\in\{\fr,\mx\}$.

\item[{(ii)}] The functions ${\mathsf p}^\ve_{w}, {\mathsf p}^\ve_{n}$ are such that
$$
{\mathsf p}^\ve_{w}\,, {\mathsf p}^\ve_{n}\,,
\sqrt{\lambda_w(x, {S}^\ve)}\, \nabla {\mathsf p}^\ve_w\,,
\sqrt{\lambda_n(x, {S}^\ve)}\, \nabla {\mathsf p}^\ve_n \in L^2(\Omega_T).
$$

\item[{(iii)}] The boundary conditions (\ref{bc3}) and the initial conditions (\ref{init1})
are satisfied.

\item[{(iv)}] For any $\varphi_w, \varphi_n \in C^1([0, T]; H^1_{\Gamma_{1}}(\Omega))$  satisfying
$\varphi_w(T) = \varphi_n(T) = 0$, we have:

\end{itemize}

\begin{small}
\begin{equation}
\label{wf-1-gl}
-\int\limits_{\Omega_{T}} \Phi^\ve(x) {S}^\ve \frac{\partial \varphi_w}{\partial t}
\, dx dt - \int\limits_{\Omega} \Phi^\ve {S}^{\bf 0} \varphi_w^{\bf 0}\, dx
+
\int\limits_{\Omega_{T}} K^\ve(x) \lambda_{w}\left(\frac{x}{\ve}, { S}^\ve\right)
\left(\nabla  {\mathsf p}^\ve_w - \vec g \right) \cdot \nabla \varphi_w\, dx dt
= \int\limits_{\Omega_{T}}  F^\ve_{w}\, \varphi_w\, dx dt;
\end{equation}
\begin{equation}
\label{gf-2-gl}
\int\limits_{\Omega_{T}} \Phi^\ve(x) {S}^\ve \frac{\partial \varphi_n}{\partial t}
\, dx dt + \int\limits_{\Omega} \Phi^\ve {S}^{\bf 0} \varphi_n^{\bf 0}\, dx
+
\int\limits_{\Omega_{T}} K^\ve(x) \lambda_{n}\left(\frac{x}{\ve}, {S}^\ve\right)
\left(\nabla  {\mathsf p}^\ve_n - \vec g \right) \cdot \nabla \varphi_n\, dx dt
= \int\limits_{\Omega_{T}}  F^\ve_{n}\,\varphi_n\, dx dt,
\end{equation}
\end{small}
where $\varphi_w^{\bf 0} \eqdef \varphi_w(0, x)$, $\varphi_n^{\bf 0} \eqdef \varphi_n(0, x)$,
and the function ${S}^{\bf 0} = {S}^{\bf 0}(x)$ is defined by the initial
condition (\ref{init1}) and the capillary pressure relation (\ref{eq.pc}).
\end{definition}

Let us also give an equivalent definition of a weak solution in terms of the global pressure and
the saturation.

\begin{definition}[\sl Weak solution in terms of global pressure and saturation]
\label{def-weak-gps}

We say that the pair of functions $\langle {S}^\ve, {\mathsf P}^\ve \rangle$ is
a weak solution of {\rm problem (\ref{gp+3-simp})} if

\begin{itemize}

\item[{(i)}] $0 \leqslant {S}^\ve \leqslant 1$ a.e. in $\Omega_T$.

\item[{(ii)}] The global pressure function
${\mathsf P}^\ve_\ell \in L^2(0, T; H^1(\Omega^\ve_\ell))$
and, for any $\ve > 0$, the saturation function $S^\ve_\ell$ is such that
$\beta_\ell(S^\ve_\ell) \in L^2(0, T; H^1(\Omega^\ve_\ell))$.

\item[{(iii)}] The boundary conditions (\ref{bc3+gp}) and the initial conditions (\ref{init1-gl+gp})
are satisfied.

\item[{(iv)}] For any $\varphi_w, \varphi_n \in C^1([0, T]; H^1_{\Gamma_{1}}(\Omega))$  satisfying
$\varphi_w(T) = \varphi_n(T) = 0$, we have:

\end{itemize}

\begin{small}
$$
-\int\limits_{\Omega_{T}} \Phi^\ve(x) {S}^\ve \frac{\partial \varphi_w}{\partial t}
\, dx dt - \int\limits_{\Omega} \Phi^\ve(x) {S}^{\bf 0} \varphi_w^{\bf 0}\, dx +
\int\limits_{\Omega^\ve_{\fr,T}} K^\ve(x)\bigg\{ \lambda_{\,\fr,w} (S^\ve_\fr)
\left(\nabla {\mathsf P}^\ve_\fr - \vec g \right) + \nabla \beta_\fr(S^\ve_\fr)\bigg\}
\cdot \nabla \varphi_w\, dx dt +
$$
\begin{equation}
\label{wf-1-gl-GP}
+ \varkappa(\ve)\,
\int\limits_{\Omega^\ve_{\mx,T}} K^\ve(x)\bigg\{ \lambda_{\,\mx,w} (S^\ve_\mx)
\left(\nabla {\mathsf P}^\ve_\mx - \vec g \right) + \nabla \beta_\mx(S^\ve_\mx)\bigg\}
\cdot \nabla \varphi_w\, dx dt
= \int\limits_{\Omega^\ve_{\fr,T}}  F^\ve_{w}\, \varphi_w\, dx dt;
\end{equation}
$$
\int\limits_{\Omega_{T}} \Phi^\ve(x) {S}^\ve \frac{\partial \varphi_n}{\partial t}
\, dx dt + \int\limits_{\Omega} \Phi^\ve(x) {S}^{\bf 0} \varphi_n^{\bf 0}\, dx +
\int\limits_{\Omega^\ve_{\fr,T}} K^\ve(x)\bigg\{ \lambda_{\,\fr,n} (S^\ve_\fr)
\left(\nabla {\mathsf P}^\ve_\fr - \vec g \right) - \nabla \beta_\fr(S^\ve_\fr)\bigg\}
\cdot \nabla \varphi_n\, dx dt +
$$
\begin{equation}
\label{gf-2-gl-GP}
+ \varkappa(\ve)\,
\int\limits_{\Omega^\ve_{\mx,T}} K^\ve(x)\bigg\{ \lambda_{\,\mx,n} (S^\ve_\mx)
\left(\nabla {\mathsf P}^\ve_\mx - \vec g \right) - \nabla \beta_\mx(S^\ve_\mx)\bigg\}
\cdot \nabla \varphi_n\, dx dt = \int\limits_{\Omega^\ve_{\fr,T}} F^\ve_{n}\, \varphi_n\, dx dt.
\end{equation}
\end{small}

\end{definition}
Existence theorem for the weak solutions defined in Definition~\ref{def-weak-simp} and
Definition~\ref{def-weak-gps} is given in \cite{APP-ex2013} in more general case of
compressible fluids.

\noindent{\bf Notational convention.} In what follows $C, C_1,..$ denote
generic constants that do not depend on $\ve$.

\section{A priori uniform estimates}

\label{uni-est}

\setcounter{equation}{0}

The uniform estimates for the initial system (\ref{debut2}) or the equivalent one
(\ref{gp+3-simp}) are given by the following lemma:

\begin{lemma}
\label{lem-uniform}
Let $\langle {\mathsf p}^\ve_w, {\mathsf p}^\ve_n , S^\ve\rangle$ be a solution to
{\rm problem (\ref{debut2})}. Then under assumptions {(A.1)-(A.9)}
the following uniform in $\ve$ estimates hold true:
$$
\big\Vert \sqrt{\lambda_{\fr,w} \left(S^\ve_\fr\right)}\, \nabla p^\ve_{\fr,w}\big\Vert_{L^2(\Omega^\ve_{\fr,T})}
+
\big\Vert \sqrt{\lambda_{\fr,n} \left(S^\ve_\fr\right)}\, \nabla p^\ve_{\fr,n} \big\Vert_{L^2(\Omega^\ve_{\fr,T})}
+
$$
\begin{equation}
\label{hhh8-cor1}
+
\varkappa^{\frac12}(\ve)\,\big\Vert \sqrt{\lambda_{\mx,w} \left(S^\ve_\mx\right)}\, \nabla p^\ve_{\mx,w}\big\Vert_{L^2(\Omega^\ve_{\mx,T})}
+
\varkappa^{\frac12}(\ve)\,\big\Vert \sqrt{\lambda_{\mx,n} \left(S^\ve_\mx\right)}\, \nabla p^\ve_{\mx,n} \big\Vert_{L^2(\Omega^\ve_{\mx,T})} \leqslant C;
\end{equation}
\begin{equation}
\label{beta-un}
\left\Vert \nabla \beta_\fr(S^\ve_\fr)\right\Vert_{L^2(\Omega^\ve_{\fr,T})}
+
\left\Vert \nabla {\mathsf P}^\ve_\fr \right\Vert_{L^2(\Omega^\ve_{\fr,T})}
+
\varkappa^{\frac12}(\ve)\,\left\Vert \nabla \beta_\mx(S^\ve_\mx)\right\Vert_{L^2(\Omega^\ve_{\mx,T})}
+
\varkappa^{\frac12}(\ve)\,\left\Vert \nabla {\mathsf P}^\ve_\mx \right\Vert_{L^2(\Omega^\ve_{\mx,T})}
\leqslant C,
\end{equation}
where $\varkappa(\ve) \eqdef \ve^\theta$ with $\theta > 0$.
\end{lemma}

\noindent{\bf Proof of Lemma \ref{lem-uniform}.}
Notice that the uniform boundedness results (\ref{hhh8-cor1}), (\ref{beta-un}) were already proved
by many authors (see, e.g., \cite{yeh2} and the references therein) in the case when the source terms
in (\ref{debut2}) were assumed to be zero. We also refer here to \cite{our-siam} and the references
therein, where the uniform boundedness results were obtained in the case of compressible two-phase
flows in porous media. Here, for reader's convenience, we recall the proof of the bounds (\ref{hhh8-cor1}),
(\ref{beta-un}) focusing on the terms involving the source functions $F^\ve_{w}, F^\ve_{n}$.

We start our analysis by obtaining the uniform bound (\ref{hhh8-cor1}). To this end
we multiply the first equation in (\ref{debut2}) by ${\mathsf p}^\ve_{w}$, the second equation in
(\ref{debut2}) by ${\mathsf p}^\ve_{n}$ and
then integrate over the domain $\Omega$. Taking into account the boundary conditions (\ref{bc3})
after integration by parts, we get the following energy equality:
$$
- \frac{d}{d t}\int\limits_\Omega \Phi^\ve(x)\,\digamma({S}^\ve) \, dx +
\int\limits_\Omega
\left\{K^\ve(x) \lambda_{w}\left(\frac{x}{\ve}, {S}^\ve\right)
\left(\nabla {\mathsf p}^\ve_{w} - \vec{g}\right)\right\} \cdot \nabla {\mathsf p}^\ve_{w} \, dx +
$$
\begin{equation}
\label{uff-2}
+
\int\limits_\Omega
\left\{K^\ve(x) \lambda_{n}\left(\frac{x}{\ve}, {S}^\ve\right)
\left(\nabla {\mathsf p}^\ve_{n} - \vec{g}\right)\right\}
\cdot \nabla {\mathsf p}^\ve_{n} \, dx
=
\int\limits_\Omega \left[F^\ve_{w}(x, t)\,{\mathsf p}^\ve_{w} + F^\ve_{n}(x, t)\,{\mathsf p}^\ve_{n}\right],
\end{equation}
where
\begin{equation}
\label{uff-2-digam}
\digamma({S}^\ve) \eqdef \digamma_\fr(S^\ve_\fr)\, {\bf 1}^\ve_\fr(x) +
\digamma_\mx(S^\ve_\mx)\, {\bf 1}^\ve_\mx(x)
\eqdef  {\bf 1}^\ve_\fr(x)\, \int\limits_1^{S^\ve_\fr} P_{\fr,c}(u)\, du +
{\bf 1}^\ve_\mx(x)\, \int\limits_1^{S^\ve_\mx} P_{\mx,c}(u)\, du.
\end{equation}

The equality (\ref{uff-2}) is the desired energy equality which will be used below to obtain the
necessary bounds that are uniform in $\ve$.  To this end we integrate (\ref{uff-2}) over the interval
$(0, T)$ to get:
$$
- \int\limits_\Omega \Phi^\ve(x)\,\digamma({S}^\ve) \, dx +
\int\limits_{\Omega_T}
\left\{K^\ve(x) \lambda_{w}\left(\frac{x}{\ve}, {S}^\ve\right)
\left(\nabla {\mathsf p}^\ve_{w} - \vec{g}\right)\right\}
\cdot \nabla {\mathsf p}^\ve_{w} \, dx dt +
$$
\begin{equation}
\label{uff-3}
+
\int\limits_{\Omega_T}
\left\{K^\ve(x) \lambda_{n}\left(\frac{x}{\ve}, {S}^\ve\right)
\left(\nabla {\mathsf p}^\ve_{n} - \vec{g}\right)\right\}
\cdot \nabla {\mathsf p}^\ve_{n} \, dx dt = \mathbb{J}^\ve_{w,n} -
\int\limits_\Omega \Phi^\ve(x)\,\digamma\left({S}^\ve(x, 0)\right) \, dx,
\end{equation}
where
\begin{equation}
\label{uff-4}
\mathbb{J}^\ve_{w,n} \eqdef \int\limits_{\Omega_T} \left[F^\ve_{w}(x, t)\,{\mathsf p}^\ve_{w} +
F^\ve_{n}(x, t)\,{\mathsf p}^\ve_{n}\right]\, dx dt.
\end{equation}

First, we notice that due to the positiveness of the porosity function $\Phi^\ve$ and the definition
of the function $\digamma\left({S}^\ve\right)$ we have that the first term on the
left-hand side of (\ref{uff-3}) is
bounded from below by a constant which does not depend on $\ve$. It is also easy to see from conditions
{(A.1)}, {(A.3)} that the second term on the right-hand side of (\ref{uff-3}) is
uniformly bounded in $\ve$. Then from (\ref{uff-3}) we get the following inequality:
$$
\int\limits_{\Omega_T}
K^\ve(x) \lambda_{w}\left(\frac{x}{\ve}, {S}^\ve\right)
\nabla {\mathsf p}^\ve_{w} \cdot \nabla {\mathsf p}^\ve_{w} \, dx dt +
\int\limits_{\Omega_T}
K^\ve(x) \lambda_{n}\left(\frac{x}{\ve}, {S}^\ve\right)
\nabla {\mathsf p}^\ve_{n} \cdot \nabla {\mathsf p}^\ve_{n} \, dx dt \leqslant
$$
\begin{equation}
\label{uff-7}
\leqslant C+
\int\limits_{\Omega_T}
K^\ve(x) \lambda_{w}\left(\frac{x}{\ve}, {S}^\ve\right) \vec{g} \cdot
\nabla {\mathsf p}^\ve_{w} \, dx dt + \int\limits_{\Omega_T}
K^\ve(x) \lambda_{n}\left(\frac{x}{\ve}, {S}^\ve\right) \vec{g} \cdot
\nabla {\mathsf p}^\ve_{n} \, dx dt + \mathbb{J}^\ve_{w,n}.
\end{equation}
With the help of Young's inequality the second and the third terms in the right-hand side of (\ref{uff-7}) can be
absorbed by the first and second term in the left-hand side of (\ref{uff-7}). Namely, we get:
\begin{equation}
\label{uff-8}
\int\limits_{\Omega_T}
K^\ve(x) \lambda_{w}\left(\frac{x}{\ve}, {S}^\ve\right)
\nabla {\mathsf p}^\ve_{w} \cdot \nabla {\mathsf p}^\ve_{w} \, dx dt +
\int\limits_{\Omega_T} K^\ve(x) \lambda_{n}\left(\frac{x}{\ve}, {S}^\ve\right)
\nabla {\mathsf p}^\ve_{n} \cdot \nabla {\mathsf p}^\ve_{n} \, dx dt
\leqslant C \big[1 + \mathbb{J}^\ve_{w,n}\big].
\end{equation}

Now it remains to estimate $\mathbb{J}^\ve_{w,n}$. Due to condition {(A.9)}, it can
be written as:
$$
\mathbb{J}^\ve_{w,n} =
\int\limits_{\Omega^\ve_{\fr,T}} \big[S^I_{\fr, w} f_I(x,t) - S^\ve_{\fr} f_P(x,t)\big]\,p^\ve_{\fr,w}
\, dx dt
+
\int\limits_{\Omega^\ve_{\fr,T}}
\big[S^I_{\fr, n} f_I(x,t) - (1 - S^\ve_{\fr}) f_P(x,t)\big] \,p^\ve_{\fr,n}\, dx dt  \eqdef
$$
\begin{equation}
\label{uff-9}
\eqdef \mathbb{J}^\ve_{w} + \mathbb{J}^\ve_{n}.
\end{equation}
Consider, first,
the term $\mathbb{J}^\ve_{w}$. From the boundedness of the saturation functions,
Cauchy's inequality and condition {(A.9)}, we get:
\begin{equation}
\label{uff-11}
\big|\mathbb{J}^\ve_{w} \big|
\leqslant \left[\Vert f_I \Vert_{L^2(\Omega_T)}
+
\Vert f_P \Vert_{L^2(\Omega_T)}\right] \Vert p^\ve_{\fr,w} \Vert_{L^2(\Omega^\ve_{\fr,T})}
\leqslant C_1\, \Vert p^\ve_{\fr,w} \Vert_{L^2(\Omega^\ve_{\fr,T})}.
\end{equation}
In a similar way,
\begin{equation}
\label{uff-12}
\big|\mathbb{J}^\ve_{n} \big| \leqslant C_2\, \Vert p^\ve_{\fr,n} \Vert_{L^2(\Omega^\ve_{\fr,T})}.
\end{equation}

Now using condition {(A.2)}, (\ref{uff-9}), (\ref{uff-11}), and (\ref{uff-12}), from the inequality
(\ref{uff-8}), we get:
$$
\mathbb{L}^\ve \eqdef
k_{\rm min}\, \int\limits_{\Omega^\ve_{\fr,T}} \lambda_{\fr,w}(S^\ve_\fr) \big|\nabla p^\ve_{\fr,w} \big|^2\, dx dt +
k_{\rm min}\, \int\limits_{\Omega^\ve_{\fr,T}} \lambda_{\fr,n}(S^\ve_\fr) \big|\nabla p^\ve_{\fr,n} \big|^2\, dx dt +
$$
$$
+
\varkappa(\ve)\, k_{\rm min}\, \int\limits_{\Omega^\ve_{\mx,T}} \lambda_{\mx,w}(S^\ve_\mx) \big|\nabla p^\ve_{\mx,w} \big|^2\, dx dt +
\varkappa(\ve)\, k_{\rm min}\, \int\limits_{\Omega^\ve_{\mx,T}} \lambda_{\mx,n}(S^\ve_\mx) \big|\nabla p^\ve_{\mx,n} \big|^2\, dx dt
\leqslant
$$
\begin{equation}
\label{uff-13}
\leqslant C_3 \left[1 + \Vert p^\ve_{\fr,w} \Vert_{L^2(\Omega^\ve_{\fr,T})} +
\Vert p^\ve_{\fr,n} \Vert_{L^2(\Omega^\ve_{\fr,T})}\right].
\end{equation}
Consider the right-hand side of (\ref{uff-13}). From (\ref{gp1}) we have:
$$
\Vert p^\ve_{\fr,w} \Vert_{L^2(\Omega^\ve_{\fr,T})} + \Vert p^\ve_{\fr,n} \Vert_{L^2(\Omega^\ve_{\fr,T})}
\leqslant
$$
\begin{equation}
\label{uff-14}
\leqslant
\left[\Vert {\mathsf P}^\ve_\fr \Vert_{L^2(\Omega^\ve_{\fr,T})} +
\Vert {\mathsf G}_{\fr,w}(S^\ve_\fr) \Vert_{L^2(\Omega^\ve_{\fr,T})} +
\Vert {\mathsf P}^\ve_\fr \Vert_{L^2(\Omega^\ve_{\fr,T})} +
\Vert {\mathsf G}_{\fr,n}(S^\ve_\fr) \Vert_{L^2(\Omega^\ve_{\fr,T})}
\right].
\end{equation}
Then, taking into account that the functions ${\mathsf G}_{\fr,w}(S^\ve_\fr), {\mathsf G}_{\fr,n}(S^\ve_\fr)$
are uniformly bounded in $\ve$, the inequality (\ref{uff-13}) takes the form:
\begin{equation}
\label{uff-15}
\mathbb{L}^\ve
\leqslant C_4 \left[1 + \Vert {\mathsf P}^\ve_\fr \Vert_{L^2(\Omega^\ve_{\fr,T})} \right].
\end{equation}
Taking into account the boundary condition ${\mathsf P}^\ve = {\mathsf P}_{\Gamma_1} = {\rm Const}$
on $\Gamma_{1} \times (0,T)$ and applying Friedrich's inequality we obtain that
\begin{equation}
\label{uff-16}
\Vert {\mathsf P}^\ve_\fr \Vert_{L^2(\Omega^\ve_{\fr,T})} \leqslant C_5\, \left[1 +
\Vert \nabla {\mathsf P}^\ve_\fr \Vert_{L^2(\Omega^\ve_{\fr,T})} \right].
\end{equation}
Finally, in view of (\ref{uff-16}), the inequality (\ref{uff-15}) takes the form:
$$
\int\limits_{\Omega^\ve_{\fr,T}} \lambda_{\fr,w}(S^\ve_\fr) \big|\nabla p^\ve_{\fr,w} \big|^2\, dx dt +
\int\limits_{\Omega^\ve_{\fr,T}} \lambda_{\fr,n}(S^\ve_\fr) \big|\nabla p^\ve_{\fr,n} \big|^2\, dx dt +
\varkappa(\ve)\, \int\limits_{\Omega^\ve_{\mx,T}} \lambda_{\mx,w}(S^\ve_\mx) \big|\nabla p^\ve_{\mx,w} \big|^2\, dx dt +
$$
\begin{equation}
\label{uff-17}
+
\varkappa(\ve)\, \int\limits_{\Omega^\ve_{\mx,T}} \lambda_{\mx,n}(S^\ve_\mx) \big|\nabla p^\ve_{\mx,n} \big|^2\, dx dt
\leqslant C_6 \left[1 + \Vert \nabla {\mathsf P}^\ve_\fr \Vert_{L^2(\Omega^\ve_{\fr,T})} \right].
\end{equation}
In order to complete the derivation of the uniform estimate, we make use of the equality
(\ref{gp6-new}). We estimate the norm
of $\nabla {\mathsf P}^\ve_\fr$ using the Cauchy inequality as follows:
\begin{equation}
\label{uff-19}
C_6\, \Vert \nabla {\mathsf P}^\ve_\fr \Vert_{L^2(\Omega^\ve_{\fr,T})} \leqslant
C_6\, \frac{\eta}{2}\, \Vert \nabla {\mathsf P}^\ve_\fr \Vert^2_{L^2(\Omega^\ve_{\fr,T})} +
C_6\, \frac{1}{2\eta},
\end{equation}
where $\eta > 0$ is an arbitrary number. Moreover, it follows from (\ref{gp6-new}) that
\begin{equation}
\label{uff-20}
\lambda_{\fr} (S^\ve_\fr) |\nabla {\mathsf P}^\ve_\fr |^2 \leqslant
\lambda_{\fr,n} (S^\ve_\fr) |\nabla p^\ve_{\fr,n}|^2 +
\lambda_{\fr,w} (S^\ve_\fr) |\nabla p^\ve_{\fr,w}|^2.
\end{equation}
Now (\ref{uff-19}) allows us to rewrite (\ref{uff-17}) in the form:
$$
\int\limits_{\Omega^\ve_{\fr,T}} \lambda_{\fr,w}(S^\ve_\fr) \big|\nabla p^\ve_{\fr,w} \big|^2\, dx dt +
\int\limits_{\Omega^\ve_{\fr,T}} \lambda_{\fr,n}(S^\ve_\fr) \big|\nabla p^\ve_{\fr,n} \big|^2\, dx dt +
\varkappa(\ve)\, \int\limits_{\Omega^\ve_{\mx,T}} \lambda_{\mx,w}(S^\ve_\mx) \big|\nabla p^\ve_{\mx,w} \big|^2\, dx dt +
$$
\begin{equation}
\label{uff-21}
+
\varkappa(\ve)\, \int\limits_{\Omega^\ve_{\mx,T}} \lambda_{\mx,n}(S^\ve_\mx) \big|\nabla p^\ve_{\mx,n} \big|^2\, dx dt
\leqslant C_6 +  C_6\, \frac{\eta}{2}\, \Vert \nabla {\mathsf P}^\ve_\fr \Vert^2_{L^2(\Omega^\ve_{\fr,T})} +
C_6\, \frac{1}{2\eta}.
\end{equation}
Let us estimate the second term on the right-hand side of (\ref{uff-21}). From condition {(A.4)}
and (\ref{uff-20}), we have:
\begin{equation}
    \begin{split}
C_6\, \frac{\eta}{2}\, \Vert \nabla {\mathsf P}^\ve_\fr \Vert^2_{L^2(\Omega^\ve_{\fr,T})}
&\leqslant
\frac{C_6\eta}{2 L_0}\, \int\limits_{\Omega^\ve_{\fr,T}} \lambda_{\fr} (S^\ve_\fr)\,
\big|\nabla {\mathsf P}^\ve_\fr \big|^2\, dx dt\\
&\leqslant
\frac{C_6\eta}{2 L_0}\, \int\limits_{\Omega^\ve_{\fr,T}}
\left[\lambda_{\fr,w}(S^\ve_\fr) \big|\nabla p^\ve_{\fr,w} \big|^2 +
\lambda_{\fr,n}(S^\ve_\fr) \big|\nabla p^\ve_{\fr,n} \big|^2 \right]\, dx dt.
    \end{split}
\label{uff-22}
\end{equation}
We set $\eta = \frac{L_0}{C_6}$ and, finally, obtain from (\ref{uff-21}) the desired inequality
(\ref{hhh8-cor1}).

Now we turn to the uniform bound (\ref{beta-un}). It immediately follows from (\ref{hhh8-cor1})
equality (\ref{gp6-new})  and the following inequality:
$\left|\nabla \beta_\ell(S^\ve_\ell) \right| \leqslant
C\,\left|\nabla \mathfrak{b}_\ell(S^\ve_\ell) \right|$.
This completes the proof of Lemma \ref{lem-uniform}. \fin

\begin{lemma}
\label{lem-uniform-ppp}
Let $\langle {\mathsf p}^\ve_w, {\mathsf p}^\ve_n, S^\ve \rangle$ be a solution to
{\rm problem (\ref{debut2})} and $\varkappa(\ve) \eqdef \ve^\theta$ with $\theta \leqslant 2$.
Then under assumptions {(A.1)-(A.9)} the following uniform in $\ve$ estimate holds true:
\begin{equation}
\label{pglob-un}
\left\Vert {\mathsf P}^\ve_\mx \right\Vert_{L^2(\Omega^\ve_{\mx,T})} \leqslant C.
\end{equation}
\end{lemma}

\noindent{\bf Proof of Lemma \ref{lem-uniform-ppp}.}
In contrast to the papers \cite{blm,yeh2},
where the standing assumptions allow to prove the continuity of the global pressure on the interface
$\Sigma^\ve_T$, in our case the global pressure is discontinuous on $\Sigma^\ve_T$. So the method
which allowed to prove (\ref{pglob-un}) by use of the extension operator from the subdomain
$\Omega^\ve_\fr$ to the whole $\Omega$ cannot be applied here. To avoid this difficulty
we make use of the ideas from \cite{ene} (see also \cite{app-2}). Since
${\mathsf P}^\ve_\mx \in L^2(0, T; H^1(\Omega^\ve_\mx))$ and
${\mathsf P}^\ve_\fr -{\mathsf P}_{\Gamma_1} \in L^2(0, T; H^1_{\Gamma_{1}}(\Omega^\ve_\fr))$, then we have:
\begin{equation}
\label{suka-1}
\Vert {\mathsf P}^\ve_\mx \Vert_{L^2(\Omega^\ve_{m,T})} \leqslant C\, \left[\ve\, \Vert \nabla {\mathsf P}^\ve_\mx \Vert_{L^2(\Omega^\ve_{\mx,T})}
+ \sqrt{\ve}\, \Vert {\mathsf P}^\ve_\mx \Vert_{L^2(\Sigma^\ve_T)} \right].
\end{equation}
Then due to the definition of the global pressure ${\mathsf P}^\ve_\mx$, (\ref{gp1}), and the interface condition
(\ref{inter-condit}) written in terms of the global pressure, one obtains the following estimate:
$$
\Vert {\mathsf P}^\ve_\mx \Vert_{L^2(\Sigma^\ve_T)} \leqslant \Vert {\mathsf P}^\ve_\mx +
{\mathsf G}_{\mx, w}(S^\ve_\mx) \Vert_{L^2(\Sigma^\ve_T)}
+
\Vert {\mathsf G}_{\mx, w}(S^\ve_\mx) \Vert_{L^2(\Sigma^\ve_T)} =
\Vert {\mathsf P}^\ve_\fr + {\mathsf G}_{\fr, w}(S^\ve_\fr) \Vert_{L^2(\Sigma^\ve_T)}
+
$$
\begin{equation}
\label{suka-2}
+ \Vert {\mathsf G}_{\mx, w}(S^\ve_\mx) \Vert_{L^2(\Sigma^\ve_T)}
\leqslant
\Vert {\mathsf P}^\ve_\fr \Vert_{L^2(\Sigma^\ve_T)}
+
\Vert {\mathsf G}_{\fr, w}(S^\ve_\fr) \Vert_{L^2(\Sigma^\ve_T)}+
\Vert {\mathsf G}_{\mx, w}(S^\ve_\mx) \Vert_{L^2(\Sigma^\ve_T)}.
\end{equation}
Now, taking into account the boundedness of ${\mathsf G}_{\ell, w}(S^\ve_\ell)$, the geometry
of $\Omega^\ve_{\mx,T}$, (\ref{suka-2}), and the estimate:
\begin{equation}
\label{suka-3}
\sqrt{\ve}\, \Vert {\mathsf P}^\ve_\fr \Vert_{L^2(\Sigma^\ve_T)} \leqslant
 C\, \left[\ve\, \Vert \nabla {\mathsf P}^\ve_\fr \Vert_{L^2(\Omega^\ve_{\fr,T})}
+ \Vert {\mathsf P}^\ve_\fr \Vert_{L^2(\Omega^\ve_{\fr,T})} \right]
\end{equation}
we obtain
\begin{equation}
\label{suka-3-1}
\Vert {\mathsf P}^\ve_\mx \Vert_{L^2(\Omega^\ve_{\mx,T})} \leqslant
C \left( \ve \Vert \nabla {\mathsf P}^\ve_\mx \Vert_{L^2(\Omega^\ve_{\mx,T})} + 1\right)
= C \left( \frac{\ve}{\varkappa^{\frac12}(\ve)} \, \varkappa^{\frac12}(\ve) \Vert \nabla {\mathsf P}^\ve_\mx \Vert_{L^2(\Omega^\ve_{\mx,T})} + 1\right).
\end{equation}
By using \eqref{beta-un}, from \eqref{suka-3-1} we get
\begin{equation}
\label{suka-3-2}
\Vert {\mathsf P}^\ve_\mx \Vert_{L^2(\Omega^\ve_{\mx,T})} \leqslant
C \left( {\ve}{\varkappa^{-\frac12}(\ve)} + 1\right),
\end{equation}
which means that for $\varkappa(\ve) \eqdef \ve^\theta$ with $\theta \leqslant 2$ the
desired inequality (\ref{pglob-un}) is obtained.
Lemma \ref{lem-uniform-ppp} is proved. \fin

Let us pass to the uniform bounds for the time derivatives of ${S}^\ve$.
In a standard way (see, e.g., \cite{ba-lp-doubpor}) we get:

\begin{lemma}
\label{cor-ps-l5}
Let $\langle {\mathsf p}^\ve_w, {\mathsf p}^\ve_n, S^\ve \rangle$ be a solution to
{\rm problem (\ref{debut2})}. Then under assumptions {(A.1)-(A.9)}
the following uniform in $\ve$ estimate holds true:
\begin{equation}
\label{theta-tt-2-cor}
\left\{\partial_t(\Phi^\ve_\ell S^\ve_\ell) \right\}_{\ve>0} \quad {\rm is \,\, uniformly \,\,
bounded \,\, in} \,\, L^2(0,T;H^{-1}(\Omega^\ve_\ell)),
\end{equation}
where the functions $\Phi^\ve_\fr, \Phi^\ve_\mx$ are defined in condition {(A.1)}.
\end{lemma}

\section{Compactness and convergence results}
\label{ex-comp-sf}

\setcounter{equation}{0}

The outline of this section is as follows. First, in subection \ref{ss-sf-1} we extend the function
$S^\ve_\fr$ from the subdomain $\Omega^\ve_\fr$ to the whole $\Omega$ and
obtain uniform estimates for the extended function $\widetilde S^\ve_\fr$. Then in
subsection \ref{ss-sf-2}, using the uniform estimates for the function $\widetilde {\mathsf P}^\ve_\fr$
 and the corresponding bounds for $\widetilde S^\ve_\fr$, we prove the compactness result for the family
$\{\widetilde S^\ve_\fr\}_{\ve>0}$. Finally, in subsection \ref{con-results}
we formulate the two-scale convergence which will be used in the derivation of the homogenized system.

\subsection{Extensions of the functions ${\mathsf P}^\ve_\fr$, $S^\ve_\fr$}
\label{ss-sf-1}

The goal of this subsection is to extend the functions ${\mathsf P}^\ve_\fr$, $S^\ve_\fr$ defined
in the subdomain $\Omega^\ve_\fr$ to the whole $\Omega$ and derive the uniform in $\ve$ estimates
for the extended functions.

Extension of the function ${\mathsf P}^\ve_\fr$.
First, we introduce the extension operator from the subdomain $\Omega^\ve_{\fr}$ to
the whole $\Omega$. Taking into account the results of \cite{acdp} we conclude that
there exists a linear continuous extension operator
$\Pi^\ve : H^1(\Omega^\ve_{\fr}) \longrightarrow H^1(\Omega)$ such that:
{(i)} $\Pi^\ve u = u$ in $\Omega^\ve_{\fr}$
and {(ii)} for any $u \in H^1(\Omega^\ve_{\fr})$,
\begin{equation}
\label{ue-15}
\Vert \Pi^\ve u \Vert_{L^2(\Omega)} \leqslant C\, \Vert u \Vert_{L^2(\Omega^\ve_{\fr})}
\quad {\rm and} \quad
\Vert \nabla (\Pi^\ve  u) \Vert_{L^2(\Omega)} \leqslant C\, \Vert \nabla u \Vert_{L^2(\Omega^\ve_{\fr})},
\end{equation}
where $C$ is a constant that does not depend on $u$ and $\ve$.
Now it follows from (\ref{beta-un}) and the Dirichlet boundary condition on $\Gamma_{1}$, that
\begin{equation}
\label{hhh9-new2}
\left\Vert \nabla (\Pi^\ve P^\ve_\fr) \right\Vert_{L^2(\Omega_{T})} +
\left\Vert \Pi^\ve P^\ve_\fr \right\Vert_{L^2(\Omega_{T})} \leqslant C.
\end{equation}

\noindent{Notational convention.}
In what follows the extension of any function $f$ will be denoted by $\widetilde f$ instead of $\Pi^\ve f$.

Extension of the function $S^\ve_\fr$.
In order to extend $S^\ve_\fr$, following
the ideas of \cite{blm}, we make use of the function $\beta_\fr$ defined in (\ref{upsi-1}).
It is evident that $\beta_\fr$ is a monotone function of $s$. Let us introduce the function:
\begin{equation}
\label{ue-23}
\beta^\ve_\fr(x, t) \eqdef \beta_\fr(S^\ve_\fr) = \int\limits_0^{S^\ve_\fr} \alpha_\fr(u)\, du.
\end{equation}
Then it follows from condition {(A.5)} that
\begin{equation}
\label{ue-24}
0 \leqslant \beta^\ve_\fr \leqslant \max_{s\in[0, 1]}
\alpha_\fr(s) \quad {\rm a.e.\,\, in} \,\, \Omega^\ve_{\fr,T}.
\end{equation}
It is also clear from (\ref{beta-un}) that
\begin{equation}
\label{ue-25}
\Vert \nabla \beta^\ve_\fr \Vert_{L^2(\Omega^\ve_{\fr,T})} \leqslant C.
\end{equation}
Hence,
\begin{equation}
\label{ue-26}
0 \leqslant \widetilde\beta^\ve_\fr \eqdef \Pi^\ve \beta^\ve_\fr \leqslant \max_{s\in[0, 1]}
\alpha_{\fr}(s) \,\, {\rm a.e.\,\, in} \,\, \Omega_T
\quad {\rm and} \quad
\Vert \nabla \widetilde \beta^\ve_\fr \Vert_{L^2(\Omega_{T})} \leqslant C.
\end{equation}

Now we can extend $S^\ve_\fr$ from $\Omega^\ve_\fr$ to the whole $\Omega$. We denote this
extension by $\widetilde S^\ve_\fr$ and define it as follows:
\begin{equation}
\label{ue-29}
\widetilde S^\ve_\fr \eqdef\, (\beta_\fr)^{-1}(\widetilde \beta^\ve_\fr).
\end{equation}
This implies that
\begin{equation}
\label{ue-31}
\int\limits_{\Omega_T} \big|\nabla \beta_\fr\big(\,\widetilde  S^\ve_\fr\,\big)\big|^2 \,\,dx\, dt =
\int\limits_{\Omega_T}
\big|\nabla \widetilde \beta^\ve_\fr \big|^2 \,\,dx\, dt \leqslant C
\quad {\rm and} \quad
0 \leqslant \widetilde  S^\ve_\fr \leqslant 1 \,\, {\rm a.e.\,\, in}
\,\, \Omega_T.
\end{equation}

\subsection{Compactness results for the sequence $\{\widetilde S^\ve_\fr\}_{\ve>0}$}
\label{ss-sf-2}

In this subsection we establish the compactness and corresponding convergence results for
the sequence $\{\widetilde  S^\ve_\fr\}_{\ve>0}$ constructed in the previous section.

\begin{proposition}
\label{prop-s}
Under our standing assumptions there is a function $S$ such that $0 \leqslant S \leqslant 1$
in $\Omega_T$ and (up to a subsequence)
\begin{equation}
\label{comp-beta-4}
\widetilde S^\ve_\fr \longrightarrow S \,\, {\rm strongly\,\, in}\,\, L^q(\Omega_T)\,\,
{\rm for\, all} \,\, 1 \leqslant q < +\infty.
\end{equation}
\end{proposition}

\noindent{\bf Proof of Proposition \ref{prop-s}.} In the proof of Proposition \ref{prop-s} we follow
the lines of \cite{blm,yeh2}. Namely, first, we establish the modulus of continuity in time
for $\widetilde \beta^\ve_\fr$ and then apply the compactness result from \cite{sim1987}.
The derivation of the modulus of continuity in time is based on the lemma obtained
earlier in \cite{yeh2}, (see also \cite{ba-lp-doubpor}).

\begin{lemma}
\label{lem-yeh-06}
For $h$ sufficiently small, we have:
\begin{equation}
\label{cor-1}
\int\limits_h^T \int\limits_{\Omega^\ve_\fr} \big[S^\ve_\fr(t) -  S^\ve_\fr(t - h)\big]\,
\big[\beta^\ve_\fr(t) - \beta^\ve_\fr(t - h) \big]\, dx\, dt \leqslant C\, h
\quad {\rm with} \,\, \beta^\ve_\fr \eqdef \beta_\fr(S^\ve_\fr),
\end{equation}
where $C$ is a constant that does not depend on $\ve, h$.

\end{lemma}

\begin{corollary}
\label{cor-yeh-06}
For $h$ sufficiently small, we have:
\begin{equation}
\label{cor-f-2}
\int\limits_{\Omega^h_T}
\big|\widetilde\beta^\ve_\fr(t) - \widetilde\beta^\ve_\fr(t - h) \big|^2\, dx\, dt
\leqslant C h
\quad  {\rm with} \,\, \Omega^h_T \eqdef \Omega \times (h, T).
\end{equation}

\end{corollary}

\noindent{\bf Proof of Corollary \ref{cor-yeh-06}.} First, let us show that
the bound (\ref{cor-1}) implies:
\begin{equation}
\label{cor-f}
\int\limits_h^T \int\limits_{\Omega^\ve_\fr}
\big|\beta^\ve_\fr(t) - \beta^\ve_\fr(t - h) \big|^2\, dx\, dt \leqslant C\, h.
\end{equation}
In fact, it is clear that due to the definition of the function $\beta_\fr$ and condition
{(A.6)} we have:
$$
\left|\beta_\fr(S^\ve_\fr(t)) - \beta_\fr(S^\ve_\fr(t-h)) \right| =
\left|\int\limits^{S^\ve_\fr(t)}_{S^\ve_\fr(t-h)} \alpha_\fr(\xi)\,d\xi \right|
\leqslant \max_{s\in[0,1]} \alpha_\fr(s)\, |S^\ve_\fr(t) -  S^\ve_\fr(t - h)|.
$$
Then from (\ref{cor-1}) we get:
$$
\int\limits_h^T \int\limits_{\Omega^\ve_\fr}
\big|\beta^\ve_\fr(t) - \beta^\ve_\fr(t - h) \big|^2\, dx\, dt \leqslant C\,
\int\limits_h^T \int\limits_{\Omega^\ve_\fr} \big[S^\ve_\fr(t) -  S^\ve_\fr(t - h)\big]\,
\big[\beta^\ve_\fr(t) - \beta^\ve_\fr(t - h) \big]\, dx\, dt \leqslant C\, h
$$
and the desired bound (\ref{cor-f}) is obtained.

Now using the property (\ref{ue-15}) of the extension operator,
from (\ref{cor-f}) we get (\ref{cor-f-2}). This completes the proof
of Corollary \ref{cor-yeh-06}. \fin

Now we are in position to complete the proof of Proposition \ref{prop-s}. First, we observe that
the sequence $\{\widetilde\beta^\ve_\fr\}_{\ve>0}$ is uniformly bounded in the space
$L^2(0, T; H^1_{\Gamma_1}(\Omega))$ and this sequence satisfies (\ref{cor-f-2}). Then it follows
from \cite{sim1987} that $\{\widetilde\beta^\ve_\fr\}_{\ve>0}$ is a compact set in the space
$L^2(\Omega_T)$ and we have that $\widetilde\beta^\ve_\fr \to \beta^\star$ strongly in
$L^2(\Omega_T)$ and due to the uniform boundedness of the function $\widetilde\beta^\ve_\fr$ in
the space $L^\infty(\Omega_T)$,
\begin{equation}
\label{comp-beta-1+}
\widetilde\beta^\ve_\fr \to \beta^\star \,\, {\rm strongly\,\, in}\,\, L^q(\Omega_T)\,\,
{\rm for\, all} \,\, 1 \leqslant q < +\infty.
\end{equation}
Now we recall that the extended saturation function $\widetilde S^\ve_\fr$ is defined
by $\widetilde S^\ve_\fr \eqdef (\beta_\fr)^{-1}(\widetilde \beta^\ve_\fr)$. We set
\begin{equation}
\label{comp-beta-3}
S \eqdef (\beta_\fr)^{-1}(\beta^\star).
\end{equation}
Then from condition {(A.6)} we have:
$$
\Vert \widetilde S^\ve_\fr - S \Vert_{L^q(\Omega_T)} =
\Vert (\beta_\fr)^{-1}(\widetilde \beta^\ve_\fr) - (\beta_\fr)^{-1}(\beta^\star) \Vert_{L^q(\Omega_T)}
\leqslant C_\beta\, \Vert \widetilde \beta^\ve_\fr - \beta^\star \Vert^{\gamma}_{L^{q\gamma}(\Omega_T)}.
$$
This inequality along with (\ref{comp-beta-1+}) implies (\ref{comp-beta-4})
and Proposition \ref{prop-s} is proved. \fin

\subsection{Two-scale convergence results}
\label{con-results}

In this subsection, taking into account the compactness results from the previous section, we formulate
the convergence results for the sequences $\{\widetilde P^\ve_\fr\}_{\ve>0}$,
$\{\widetilde S^\ve_\fr\}_{\ve>0}$. In this paper the homogenization process for the problem is
rigorously obtained by using the two-scale approach, see, e.g., \cite{al}.
For the reader's convenience, let us recall the definition of the two-scale convergence.

\begin{definition}
\label{two}
A sequence of functions $\{v^\ve\}_{\ve>0} \subset L^2(\Omega_T)$ two-scale converges to
$v \in L^2(\Omega_T\times Y)$ if $\Vert v^\ve\Vert_{L^2(\Omega_T)} \leqslant C$,
and for any test function $\varphi \in C^\infty(\overline{\Omega_T}; C_\#(Y))$ the following relation holds:
$$
\lim_{\ve\to 0} \int\limits_{\Omega_T} v^\ve(x, t)\, \varphi \left(x, \frac{x}{\ve}, t\right)\, dx\, dt =
\int\limits_{\Omega_T \times Y} v(x, y, t)\,
\varphi(x, y, t) \, dy\, dx\, dt.
$$
\end{definition}

\noindent
This convergence is denoted by $v^\ve(x, t) \stackrel {2s}
\rightharpoonup v(x, y, t)$.

\medskip

Following \cite{al-95} we also introduce {\it the two-scale convergence on periodic surfaces:}
\begin{definition}
\label{two-surf}
A sequence of functions $\{v^\ve\}_{\ve>0} \subset L^2(\Sigma^{\ve}_T)$
two-scale converges to $v \in L^2(\Omega_T; L^2(\Gamma_{\fr\mx}))$
on $\Gamma_{\fr\mx}$ if  for any test function
$\varphi \in C^\infty(\overline{\Omega_T}; C_\#(Y))$ the following relation holds:
$$
\lim_{\ve\to 0} \ve \int\limits_{\Sigma^{\ve}_T} v^\ve(x, t)\,
\varphi \left(x, \frac{x}{\ve}, t\right)\, dH^{d-1}(x)\, dt =
\int\limits_{\Omega_T}\int\limits_{\Gamma_{\fr\mx}} v(x, y, t)\,
\varphi(x, y, t) \, dH^{d-1}(y)\, dx\, dt,
$$
where, as before $\Sigma^{\ve}_T \eqdef \Gamma_{\fr\mx}^{\ve}\times (0,T)$, and $dH^{d-1}$ is
the $(d-1)$-dimensional Hausdorff measure.
\end{definition}
\noindent
This convergence is denoted by $v^\ve(x, t) \stackrel {2s-\Gamma_{\mx\fr}}
\rightharpoonup v(x, y, t)$.
\medskip

Now we summarize the convergence results for the sequences $\{\widetilde P^\ve_\fr\}_{\ve>0}$
and $\{\widetilde S^\ve_\fr\}_{\ve>0}$. We have:

\begin{lemma}
\label{2scale}
For any rate of contrast there exist a function $S$ such that $0 \leqslant S \leqslant 1$
a.e. in $\Omega_T$,
$\beta_\fr(S)-\beta_\fr(1) \in L^2(0, T; H^1_{\Gamma_1}(\Omega))$,
and functions
${\mathsf P}-{\mathsf P}_{\Gamma_1} \in L^2(0, T; H^1_{\Gamma_1}(\Omega))$,
${\mathsf w}_p, {\mathsf w}_s \in L^2(\Omega_T; H^1_{per}(Y))$ such that up to a subsequence:
\begin{equation}
\label{2s-1}
\widetilde S^\ve_\fr(x, t) \longrightarrow S(x, t) \,\,
{\rm strongly\,\, in}\,\, L^q(\Omega_T)\,\, \forall \ 1 \leqslant q < +\infty;
\end{equation}
\begin{equation}
\label{2s-2}
\widetilde {\mathsf P}^\ve_\fr(x, t) \rightharpoonup {\mathsf P}(x, t) \,\,
{\rm weakly\,\, in}\,\, L^2(0, T; H^1(\Omega));
\end{equation}
\begin{equation}
\label{2s-3}
\nabla \widetilde {\mathsf P}^\ve_\fr(x, t)\stackrel {2s}
\rightharpoonup \nabla {\mathsf P}(x, t) + \nabla_y {\mathsf w}_p(x, t, y);
\end{equation}
\begin{equation}
\label{2s-40}
\beta_\fr(\widetilde S^\ve_\fr) \longrightarrow \beta_\fr(S) \,\,
{\rm strongly\,\, in}\,\, L^q(\Omega_T)\,\, \forall \ 1 \leqslant q < +\infty;
\end{equation}
\begin{equation}
\label{2s-4}
\nabla \beta_\fr(\widetilde S^\ve_\fr) (x, t) \stackrel {2s} \rightharpoonup
\nabla \beta_\fr(S)(x, t) + \nabla_y {\mathsf w}_s(x, t, y);
\end{equation}
\begin{equation}
\label{2s-5}
\widetilde {\mathsf P}^\ve_\fr(x, t)  \stackrel {2s-\Gamma_{\mx\fr}}\rightharpoonup {\mathsf P}(x, t);
\end{equation}
\begin{equation}
\label{2s-6}
\beta_\fr(\widetilde S^\ve_\fr(x, t))  \stackrel {2s-\Gamma_{\mx\fr}}\rightharpoonup \beta_\fr(S(x, t)) .
\end{equation}

\end{lemma}

\noindent
The {\bf Proof of Lemma~\ref{2scale}} is based on the {\it a priori estimates}
for the functions $\beta_\fr(S^\ve_\fr)$ and ${\mathsf P}^\ve_\fr$ obtained in Section \ref{uni-est},
the extension results from Subsection \ref{ss-sf-1}, and Proposition \ref{prop-s}.
The two-scale convergence results (\ref{2s-3}) and (\ref{2s-4}) are
obtained by arguments similar to those in \cite{al}. The two-scale convergence (\ref{2s-5})
and (\ref{2s-6}) can be proved by applying Proposition 2.6 in \cite{al-95}.
Lemma \ref{2scale} is proved. \fin

Note also that the notion of strong two-scale convergence on periodic surfaces can be introduced
in analogy with the ordinary strong two-scale convergence.

\begin{definition}
\label{str-conv-per-surf}
A sequence $\{v^\ve\}_{\ve>0} \subset L^2(\Sigma^{\ve}_T)$ converges the two-scale strongly
to $v \in L^2(\Omega_T; L^2(\Gamma_{\fr\mx}))$ on $\Gamma_{\fr\mx}$  if
$$
\lim_{\ve\to 0} \ve \int\limits_{\Sigma^{\ve}_T} | v^\ve(x, t)- v\left(x, \frac{x}{\ve}, t\right)|^2\,
dH^{d-1}(x)\, dt = 0.
$$
\end{definition}

It is easy to verify that the strong two-scale convergence on periodic surfaces implies
the two-scale convergence on periodic surfaces with the same limit.

Using the strong convergence (\ref{2s-40}) and the boundedness of
$\nabla\beta_\fr(\widetilde S^\ve_\fr)$ given in Lemma~\ref{lem-uniform} we get:
\begin{align*}
\ve \| \beta_\fr(\widetilde S^\ve_\fr) - \beta_\fr(S)\|_{L^2(\Sigma^{\ve}_T)}^2 \leqslant C \left[
\ve^2 \| \nabla  \beta_\fr(\widetilde S^\ve_\fr)  -\nabla \beta_\fr(S)\|_{L^2(\Omega_{\fr,T}^\ve)}^2 +
\|\beta_\fr(\widetilde S^\ve_\fr) - \beta_\fr(S) \|_{L^2(\Omega_{\fr,T}^\ve)}^2
\right],
\end{align*}
which tends to zero on a given subsequence as $\ve\to 0$. Therefore, we conclude that the sequence
$\{\beta_\fr(\widetilde S^\ve_\fr)\}_{\ve>0}$
converges strongly two-scale on the surface $\Gamma_{\fr\mx}$ to $\beta_\fr(S)$. Furthermore, we have:
\begin{lemma}
Let  $\{\beta_\fr(\widetilde S^\ve_\fr)\}$ be a subsequence from Lemma~\ref{2scale}. Then for any
Lipschitz function ${\cal M}\colon [0, \beta_\fr(1)]\to \mathbb{R}$
the sequence  $\{{\cal M}(\beta_\fr(\widetilde S^\ve_\fr))\}_{\ve>0}$ converges strongly
two-scale on the surface $\Gamma_{\fr\mx}$ to ${\cal M}(\beta_\fr(S))$.
\label{lema:2sc-surf}
\end{lemma}
Lemma~\ref{lema:2sc-surf} follows immediately from the estimate
\begin{align*}
     \| {\cal M}(\beta_\fr(\widetilde S^\ve_\fr)) - {\cal M}(\beta_\fr(S))\|_{L^2(\Sigma^{\ve}_T)}^2 \leq
     L_{\cal M}^2 \| \beta_\fr(\widetilde S^\ve_\fr) - \beta_\fr(S)\|_{L^2(\Sigma^{\ve}_T)}^2,
\end{align*}
where $L_{\cal M}$ is the Lipschitz constant which does not depend on $\ve$.

\section{Dilation operator and convergence results}
\label{dil-oper}

\setcounter{equation}{0}

It is known that due to the nonlinearities and the strong coupling of the problem,
the two-scale convergence does not provide an explicit form for the source terms
appearing in the homogenized model, see for instance \cite{blm,choq,yeh2}.
To overcome this difficulty the authors make use of the dilation operator.
Here we refer to \cite{adh,blm,choq,yeh2} for the definition and main properties
of the dilation operator. Let us also notice that the notion of the dilation
operator is closely related to the notion of the unfolding operator. We refer here,
e.g., to \cite{ddg} for the definition and the properties of this operator.

The outline of this section is as follows. First, in subsection \ref{def-base-dil} we introduce
the definition of the dilation operator and describe its main properties.
Then in subsection \ref{dil-func-prop} we obtain the equations for the dilated saturation and
the global pressure functions and the corresponding uniform estimates.
Finally, in subsection \ref{dil-func-conv} we consider the convergence results for the dilated functions.

\subsection{Definition and preliminary results}
\label{def-base-dil}

\begin{definition}
\label{def-dilop}
For a given $\ve > 0$, we define a dilation operator $\mathfrak{D}^\ve$
mapping measurable functions defined in $\Omega^\ve_{\mx,T}$ to measurable
functions defined in $\Omega_T \times Y_\mx$ by
\begin{equation}
\label{do-1}
\left(\mathfrak{D}^\ve \varphi \right)(x, y, t) \eqdef
\left\{
\begin{array}[c]{ll}
\varphi\left( c^\ve(x) + \ve\, y, t\right),
\quad {\rm if}\,\, c^\ve(x) + \ve\, y \in \Omega^\ve_\mx; \\[2mm]
0,
\quad {\rm elsewhere}, \\
\end{array}
\right.
\end{equation}
where $c^\ve(x) \eqdef \ve\, k$ if $x \in \ve\, (Y + k)$
with $k \in \mathbb{Z}^d$ denotes the lattice translation point of the $\ve$-cell
domain containing $x$.
\end{definition}

The basic properties of the dilation operator are given by the following lemma (see \cite{adh,yeh2}).

\begin{lemma}
\label{lemma-dilop1}

Let $\varphi, \psi \in L^2(0, T; H^1(\Omega^\ve_\mx))$. Then we have:
\begin{equation}
\label{dilll}
\nabla_y \mathfrak{D}^\ve \varphi = \ve\, \mathfrak{D}^\ve (\nabla_x \varphi)
\quad {\rm a.e.\,\,in}\,\,\Omega_T \times Y_\mx;
\end{equation}
$$
\Vert \mathfrak{D}^\ve \varphi \Vert_{L^2(\Omega_T \times Y_\mx)} =
\Vert \varphi \Vert_{L^2(\Omega^\ve_{\mx,T})}; \,\,
\Vert \nabla_y \mathfrak{D}^\ve \varphi \Vert_{L^2(\Omega_T \times Y_\mx)}
= \ve\,
\Vert \mathfrak{D}^\ve \nabla_x\, \varphi \Vert_{L^2(\Omega_T \times Y_\mx)}
= \ve\,
\Vert \nabla_x\, \varphi \Vert_{L^2(\Omega^\ve_{\mx,T})};
$$
$$
\left(\mathfrak{D}^\ve \varphi, \mathfrak{D}^\ve \psi
\right)_{L^2(\Omega_T \times Y_\mx)} =
\left(\varphi, \psi \right)_{L^2(\Omega^\ve_{\mx,T})}.
$$
\end{lemma}

The following lemma gives the link between the two-scale and the weak convergence
(see, e.g., \cite{blm}).

\begin{lemma}
\label{lemma-dilop}
Let $\{\varphi^\ve\}_{\ve>0}$ be a uniformly bounded sequence in $L^2(\Omega^\ve_{\mx,T})$
satisfying: {(i)} $\mathfrak{D}^\ve \varphi^\ve \rightharpoonup \varphi^0$
weakly in $L^2(\Omega_T; L^2_{per}(Y_\mx))$; {(ii)} ${\bf 1}^\ve_\mx(x)
\varphi^\ve \stackrel {2s} \rightharpoonup \varphi^* \in L^2(\Omega_T; L^2_{per}(Y_\mx))$.
Then $\varphi^0 = \varphi^*$ a.e. in $\Omega_T \times Y_\mx$.

\end{lemma}

Finally, we also have the following result (see, e.g., \cite{choq,yeh2}).

\begin{lemma}
\label{lemma-dilop2}
If $\varphi^\ve \in L^2(\Omega^\ve_{\mx,T})$ and
${\bf 1}^\ve_\mx(x) \varphi^\ve \stackrel {2s} \to \varphi \in L^2(\Omega_T; L^2_{per}(Y_\mx))$
then $\mathfrak{D}^\ve \varphi^\ve$ converges to $\varphi$ strongly in $L^2(\Omega_T \times Y_\mx)$.
Here $\stackrel {2s} \to$ denotes the strong two-scale convergence.
If $\varphi \in L^2(\Omega_T)$ is considered as an element of $L^2(\Omega_T \times Y_\mx)$
constant in $y$, then $\mathfrak{D}^\ve \varphi$ converges strongly to $\varphi$ in
$L^2(\Omega_T \times Y_\mx)$.
\end{lemma}

The dilation operator shows the same properties with respect to the two-scale convergence on
periodic surfaces. For a given function $v\in  L^2(\Sigma^{\ve}_T)$ and from definition
of the dilation operator we have $\mathfrak{D}^\ve(v)\in L^2(\Omega_T; L^2(\Gamma_{\fr\mx}))$
and
$$
\sqrt{\ve}\| v\|_{L^2(\Sigma^{\ve}_T)} = \| \mathfrak{D}^\ve(v)\|_{L^2(\Omega_T; L^2(\Gamma_{\fr\mx}))}.
$$
We have also the following lemma:

\begin{lemma}
If $\{v^\ve\}_{\ve>0} \subset L^2(\Sigma^{\ve}_T)$ is a sequence that converges to
$v \in L^2(\Omega_T; L^2(\Gamma_{\fr\mx}))$ in the two-scale sense on $\Gamma_{\fr\mx}$,
then the sequence $\{\mathfrak{D}^\ve(v^\ve)\}_{\ve>0}$ converges weakly to the same limit, that is
$\mathfrak{D}^\ve(v^\ve) \rightharpoonup v$ in $L^2(\Omega_T; L^2(\Gamma_{\fr\mx}))$.
If $\{v^\ve\}_{\ve>0} \subset L^2(\Sigma^{\ve}_T)$  converges strongly  to
$v \in L^2(\Omega_T; L^2(\Gamma_{\fr\mx}))$  in the two-scale sense on $\Gamma_{\fr\mx}$,
then the sequence $\{\mathfrak{D}^\ve(v^\ve)\}_{\ve>0}$ converges strongly to the same limit in
$L^2(\Omega_T; L^2(\Gamma_{\fr\mx}))$.
\label{lemma-dilop-3}
\end{lemma}

Due to Lemma~\ref{lema:2sc-surf}, one can apply  Lemma~\ref{lemma-dilop-3} to the  sequence
$\{ {\cal M}(\beta_\fr(\widetilde S^\ve_\fr))\}_{\ve>0}$ and find a subsequence, such that
\begin{align*}
\int\limits_{\Omega_T} \int\limits_{\Gamma_{\fr\mx}} \left|
{\cal M}(\beta_\fr(\mathfrak{D}^\ve(\widetilde S^\ve_\fr))) - {\cal M}(\beta_\fr(S))
\right|^2 \, dH^{d-1}(y)\, dx\, dt \to 0
\end{align*}
when $\ve\to 0$, for any Lipschitz function ${\cal M}$. As a consequence
we have.

\begin{corollary}
\label{corol-dilop}
Let ${\cal M}\colon [0, \beta_\fr(1)]\to \mathbb{R}$
be a Lipschitz function.
Then there is a subsequence $\ve = \ve_k$ of the sequence
$\{{\cal M}(\beta_\fr(\widetilde S^\ve_\fr))\}_{\ve>0}$, still denoted by $\ve$,
such that for a.e. $x\in\Omega$
\begin{align*}
\int\limits_0^T \int\limits_{\Gamma_{\fr\mx}} \left|
{\cal M}(\beta_\fr(\mathfrak{D}^\ve(\widetilde S^\ve_\fr(x,y,t)))) - {\cal M}(\beta_\fr(S(x,y,t)))
\right|^2 \, dH^{d-1}(y) dt \to 0  \quad\text{as }\;  \ve\to 0.
\end{align*}
\end{corollary}

\subsection{The dilated functions $\mathfrak{D}^\ve S^\ve_\mx, \mathfrak{D}^\ve P^\ve_\mx$
and their properties}
\label{dil-func-prop}

In this section we derive the equations for the dilated functions $\mathfrak{D}^\ve S^\ve_\mx,
\mathfrak{D}^\ve P^\ve_\mx$ and obtain the corresponding uniform estimates.
In what follows we also make use of the notation:
$$
\mathfrak{D}^\ve S^\ve_\mx \eqdef s^\ve_\mx \quad {\rm and} \quad
\mathfrak{D}^\ve {\mathsf P}^\ve_\mx \eqdef p^\ve_\mx.
$$
The equations for the dilated functions $s^\ve_\mx, p^\ve_\mx$ are given by the following lemma.

\begin{lemma}
\label{lem-dil}
For $x \in \Omega$, the functions $s^\ve_\mx, p^\ve_\mx$
satisfy the following system of equations:
\begin{equation}
\label{dilo-1}
\Phi_\mx(y) \frac{\partial s^\ve_\mx}{\partial t} - \frac{\varkappa(\ve)}{\ve^2}\,
{\rm div}_y\, \bigg\{K(x, y) \left[
\lambda_{\mx,w} (s^\ve_\mx) \nabla_y p^\ve_\mx +
\nabla_y \beta_\mx(s^\ve_\mx) - \ve\,\lambda_{\mx,w} (s^\ve_\mx) \vec{g}\, \right] \bigg\}
= 0;
\end{equation}
\begin{equation}
\label{dilo-2}
- \Phi_\mx(y) \frac{\partial s^\ve_\mx}{\partial t}
- \frac{\varkappa(\ve)}{\ve^2}\, {\rm div}_y\, \bigg\{K(x, y)\, \left[ \lambda_{\mx,n} (s^\ve_\mx)
\nabla_y p^\ve_\mx - \nabla_y \beta_\mx(s^\ve_\mx)
- \ve\,\lambda_{\mx,n}(s^\ve_\mx)\,\vec{g}\,
\right] \bigg\} = 0,
\end{equation}
in the space $L^2(0, T; H^{-1}(Y_\mx))$.

\end{lemma}

\noindent The {\bf Proof of Lemma \ref{lem-dil}} is given in \cite{blm,yeh2}.

The system of equations (\ref{dilo-1})-(\ref{dilo-2}) is provided with the following
boundary conditions:
\begin{equation}
\label{dilo-6}
\beta_\mx(s^\ve_\mx) = {\cal M}(\beta_\fr(\mathfrak{D}^\ve \tilde{S^\ve_\fr}))
\quad {\rm on}\; \Gamma_{\fr\mx}
\end{equation}
for $(x, t) \in \Omega^\ve_\mx \times (0,T)$, where
\begin{equation}
{\cal M} \eqdef \beta_{\mx} \circ ( P_{\mx,c})^{-1}\circ P_{\fr,c} \circ (\beta_{\fr})^{-1}.
\label{function:M}
\end{equation}
Note that under our hypothesis function $ {\cal M}$ is Lipschitz continuous.
We also have
\begin{equation}
\label{dilo-6a}
{\mathsf p}^\ve_\mx + {\mathsf G}_{\mx, w}(s^\ve_\mx) =
\mathfrak{D}^\ve {\mathsf P}^\ve_\fr + {\mathsf G}_{\fr, w}(\mathfrak{D}^\ve \tilde{S^\ve_\fr})
\,\,
{\rm and} \,\,
{\mathsf p}^\ve_\mx + {\mathsf G}_{\mx,n}(s^\ve_\mx)
= \mathfrak{D}^\ve{\mathsf P}^\ve_\fr + {\mathsf G}_{\fr,n}(\mathfrak{D}^\ve\tilde{S^\ve_\fr})
\end{equation}
on $\Gamma_{\fr\mx}$ for $(x, t) \in \Omega^\ve_\mx \times (0,T)$.

The initial conditions are
\begin{equation}
\label{dilo-7}
s^\ve_\mx(x, y, 0) = (\mathfrak{D}^\ve S^{\bf 0}_\mx)(x, y) \quad {\rm and} \quad
p^\ve_\mx(x, y, 0) = (\mathfrak{D}^\ve {\mathsf P}^{\bf 0}_\mx)(x, y) \quad {\rm in} \,\,
\Omega^\ve_\mx \times Y_\mx,
\end{equation}
where $S^{\bf 0}_\mx, {\mathsf P}^{\bf 0}_\mx$ are the restrictions to the domain $\Omega^\ve_\mx$
of the functions $S^{\bf 0}, {\mathsf P}^{\bf 0}$  defined in (\ref{init1-gl+gp}) and
the dilations of the functions defined on the fracture system can be defined in a way
similar to one already used for the functions defined on the matrix part.

Now we establish {\it a priori} estimates for the functions $s^\ve_\mx, p^\ve_\mx$.
They are given by the following lemma.

\begin{lemma}
\label{lem-dil-est}
Let $\langle s^\ve_\mx, p^\ve_\mx \rangle$ be a solution to {\rm problem
(\ref{dilo-1})-(\ref{dilo-2})}. Then:
\begin{itemize}

\item[{(i)}] {For any rate of contrast ($\theta > 0$)},
\begin{equation}
\label{dilo-8}
0 \leqslant s^\ve_\mx \leqslant 1 \quad {\rm a.e.\,\, in\,\,} \Omega_T \times Y_\mx;
\end{equation}
\begin{equation}
\label{dilo-12}
\Vert \partial_t(\Phi_\mx\, s^\ve_\mx) \Vert_{L^2(\Omega_T; H^{-1}_{per}(Y_\mx))}
\leqslant C.
\end{equation}

\item[{(ii)}] {For the high contrast in the critical case ($\theta = 2$)},
\begin{equation}
\label{dilo-9}
\Vert \nabla_y \beta_\mx(s^\ve_\mx) \Vert_{L^2(\Omega_T; L^2_{per}(Y_\mx))} \leqslant C;
\end{equation}
\begin{equation}
\label{dilo-10}
\Vert p^\ve_\mx \Vert_{L^2(\Omega_T; H^1_{per}(Y_\mx))} \leqslant C.
\end{equation}

\item[{(iii)}] {For the moderate contrast ($0 < \theta < 2$)},
\begin{equation}
\label{dilo-10<2}
\ve^{\frac{\theta}{2}-1}\,
\Vert \nabla_y \beta_\mx(s^\ve_\mx) \Vert_{L^2(\Omega_T; L^2_{per}(Y_\mx))}
+
\Vert \beta_\mx(s^\ve_\mx) \Vert_{L^2(\Omega_T; L^2_{per}(Y_\mx))}
\leqslant C;
\end{equation}
\begin{equation}
\label{dilo-10<2p}
\ve^{\frac{\theta}{2}-1}\,
\Vert \nabla_y p^\ve_\mx \Vert_{L^2(\Omega_T; L^2_{per}(Y_\mx))}
+
\Vert p^\ve_\mx \Vert_{L^2(\Omega_T; L^2_{per}(Y_\mx))}
\leqslant C.
\end{equation}
\end{itemize}
\end{lemma}

\noindent{\bf Proof of Lemma \ref{lem-dil-est}.} Statement (\ref{dilo-8}) is evident.
The bound (\ref{dilo-12}) with $\Phi_\mx = \Phi_\mx(y)$ follow from Lemma \ref{cor-ps-l5} and
Lemma \ref{lemma-dilop1}. The uniform estimates for $p^\ve_\mx$ in (\ref{dilo-10}) and
(\ref{dilo-10<2p}) follow from the uniform bound (\ref{pglob-un}) and Lemma \ref{lemma-dilop1}.
The uniform estimates for the gradients of the functions $\beta_\mx(s^\ve_\mx)$ and
$p^\ve_\mx$ easy follow from the uniform bounds (\ref{beta-un}) and Lemma \ref{lemma-dilop1}.
Lemma \ref{lem-dil-est} is proved. \fin

\begin{remark}
\label{rem-supercrit-1}
Notice that in what follows we do not need the uniform estimates for the
dilated functions in the case of the very high contrast.
\end{remark}

\subsection{Convergence results for the dilated functions}
\label{dil-func-conv}

In this subsection we establish convergence results which will be used below to obtain
the homogenized system. From Lemmas \ref{lemma-dilop}, \ref{lem-dil-est} we get
the following convergence results.

\begin{lemma}
\label{conv-lemma-dil}
Let $\langle s^\ve_\mx, p^\ve_\mx \rangle$ be a solution to {\rm problem
(\ref{dilo-1})-(\ref{dilo-2}), (\ref{dilo-6})-(\ref{dilo-7})}. Then (up to a subsequence),
\begin{itemize}

\item[{(i)}] {For the high contrast in the critical case ($\theta = 2$)},
\begin{equation}
\label{weak-dil-10}
{\bf 1}^\ve_\mx(x) S^\ve_\mx \stackrel {2s} \rightharpoonup s \in L^2(\Omega_T; L^2_{per}(Y_\mx))
\quad {\rm and} \quad s^\ve_\mx \rightharpoonup s \,\, {\rm weakly\,\, in}\,\, L^2(\Omega_T \times Y_\mx);
\end{equation}
\begin{equation}
\label{weak-dil-2}
{\bf 1}^\ve_\mx(x)\, {\mathsf P}^\ve_\mx \stackrel {2s} \rightharpoonup p \in L^2(\Omega_T; L^2_{per}(Y_\mx))
\quad {\rm and} \quad
p^\ve_\mx \rightharpoonup p \,\, {\rm weakly\,\, in}\,\, L^2(\Omega_T; H^1(Y_\mx));
\end{equation}
\begin{equation}
\label{weak-dil-201}
{\bf 1}^\ve_\mx(x)\, \nabla_x {\mathsf P}^\ve_\mx \stackrel {2s} \rightharpoonup \nabla_y p \in L^2(\Omega_T; L^2_{per}(Y_\mx));
\end{equation}
\begin{equation}
\label{weak-dil-30}
{\bf 1}^\ve_\mx(x) \beta_\mx(S^\ve_\mx) \stackrel {2s} \rightharpoonup \beta^*
\quad {\rm and} \quad
\beta_\mx(s^\ve_\mx) \rightharpoonup \beta^* \,\, {\rm weakly\,\, in}\,\, L^2(\Omega_T; H^1(Y_\mx));
\end{equation}
\begin{equation}
\label{weak-dil-301}
{\bf 1}^\ve_\mx(x)\, \nabla_x \beta_\mx(S^\ve_\mx) \stackrel {2s} \rightharpoonup
\nabla_y \beta^*.
\end{equation}

\item[{(ii)}] {For the very high contrast ($\theta > 2$)},
\begin{equation}
\label{weak-dil-10very}
{\bf 1}^\ve_\mx(x) S^\ve_\mx \stackrel {2s} \rightharpoonup s \in L^2(\Omega_T; L^2_{per}(Y_\mx)).
\end{equation}

\item[{(iii)}] {For the moderate contrast ($0 < \theta < 2$)},
\begin{equation}
\label{weak-dil-10mod}
{\bf 1}^\ve_\mx(x) S^\ve_\mx \stackrel {2s} \rightharpoonup s \in L^2(\Omega_T; L^2_{per}(Y_\mx))
\quad {\rm and} \quad s^\ve_\mx \rightharpoonup s \,\, {\rm weakly\,\, in}\,\, L^2(\Omega_T \times Y_\mx);
\end{equation}
\begin{equation}
\label{weak-dil-30mod}
{\bf 1}^\ve_\mx(x) \beta_\mx(S^\ve_\mx) \stackrel {2s} \rightharpoonup \beta^*_1
\quad {\rm and} \quad
{\bf 1}^\ve_\mx(x)\, \ve^{\theta} \nabla \beta_\mx(S^\ve_\mx) \stackrel {2s}
\rightharpoonup \beta_1;
\end{equation}
\begin{equation}
\label{weak-dil-301mod}
\beta_\mx(s^\ve_\mx) \rightharpoonup \beta^*_1 \,\, {\rm weakly\,\, in}\,\, L^2(\Omega_T; H^1(Y_\mx)).
\end{equation}

\end{itemize}
\end{lemma}

It is important to notice that the convergence results of
Lemma \ref{conv-lemma-dil} are not sufficient for derivation of the equations for the
limit functions $\langle s, p \rangle$ which involve only these functions and not the undefined
limits appearing in (\ref{weak-dil-30}), (\ref{weak-dil-301}), (\ref{weak-dil-30mod})
and (\ref{weak-dil-301mod}). In order to overcome this difficulty,
we introduce the restrictions of the functions $s^\ve_\mx$, $p^\ve_\mx$ which are defined below.
For these functions we obtain more estimates which allow us to obtain the desired equations.
For this, we make use of the density arguments. Namely, following \cite{choq}
(see also \cite{ba-lp-doubpor}), we fix $x \in \Omega$ and define the restrictions of
$s^\ve_\mx$, $p^\ve_\mx$ to the $\ve$-cell containing the point $x$.
These functions are defined in the domain $Y_\mx \times (0, T)$ and are
constants in the slow variable $x$. In order to obtain the uniform estimates for the
restricted functions (they are similar to the corresponding estimates for
${\mathsf P}^\ve_\fr$, $S^\ve_\fr$ from Section \ref{uni-est}) we make use of the estimates
(\ref{dilo-8})-(\ref{dilo-10<2p}).

The scheme is as follows. First, for any natural ${\bf n}$, we introduce the set of points
$x \in \Omega$ such that the corresponding norms for the restricted functions are not uniformly bounded
in $\ve$. It turns out that the measure of this set is asymptotically small as ${\bf n} \to +\infty$
(see Propositions \ref{proppi-1}, \ref{proppi-2} below). Then taking into account this fact and
using the estimates (\ref{dilo-8})-(\ref{dilo-10<2p}), we, finally, obtain the desired uniform
estimates for the restricted functions  (see Lemma \ref{lem-hitro-vyeb} below).

Let us first denote a periodicity cell $ \ve \big(Y +  k\big)$  which contains point
$x_0$ by $K^\ve_{x_0}$. For given $x_0$ and $\ve$ the index $k\in \mathbb{Z}^d$ which defines the cell $K^\ve_{x_0}$
can be uniquely defined and therefore we have a well defined function  $k(x_0, \ve) \in \mathbb{Z}^d$ such that
$K^\ve_{x_0} \eqdef \ve \big(Y +  k(x_0, \ve)\big)$.
Due to the definition of the dilation operator dilated  functions are constant in $x$ on $K^\ve_{x_0}$.
The restricted functions are given by:
\begin{equation}
\label{new-smx}
s^\ve_{\mx,x_0}(y, t) \eqdef
\left\{
\begin{array}[c]{ll}
s^\ve_\mx,
\quad {\rm for}\,\, x \in K^\ve_{x_0}; \\[2mm]
0, \quad {\rm if \,\, not}; \\
\end{array}
\right.
\quad
p^\ve_{\mx,x_0}(y, t) \eqdef
\left\{
\begin{array}[c]{ll}
p^\ve_\mx,
\quad {\rm for}\,\, x \in K^\ve_{x_0}; \\[2mm]
0, \quad {\rm if \,\, not}. \\
\end{array}
\right.
\end{equation}
For any $\ve > 0$, the pair $\langle s^\ve_{\mx,x_0}, p^\ve_{\mx,x_0} \rangle$ is a solution to
problem (\ref{dilo-1})-(\ref{dilo-2}), (\ref{dilo-6})-(\ref{dilo-7}) in $Y_\mx \times (0,T)$.

Now we estimate the measure of the set of points $x \in \Omega$ such that the corresponding norms
for the restricted functions are not uniformly bounded in $\ve$. The following result holds true.

\begin{proposition}
\label{proppi-1}
Let $f^\ve_\mx = f^\ve_\mx(x, y, t)$ be a dilated function such that
\begin{equation}
\label{rev-fff}
\Vert f^\ve_\mx \Vert_{L^2(\Omega_T; L^2_{per}(Y_\mx))} \leqslant C
\end{equation}
and let $A_{\bf n}$ be a set of points defined by
\begin{equation}
\label{rev-fff-5}
A_{\bf n} \eqdef \left\{x \in \Omega\,:\, \liminf_{\ve\to0}\Vert \widehat f^\ve_{\mx,k(x,\ve)}
\Vert_{L^2(0,T; L^2_{per}(Y_\mx))} \geqslant {\bf n} \right\},
\end{equation}
where for fixed $k\in \mathbb{Z}^d$
\begin{equation}
\label{rev-fff-2}
\widehat f^\ve_{\mx,k}(y, t) \eqdef
\left\{
\begin{array}[c]{ll}
f^\ve_\mx(\ve k, y, t),
\quad {\rm if}\,\, k \,\, {\rm is \,\, such\,\, that}\,\,\,
\ve(Y_\mx + k) \cap \Omega \not=\emptyset; \\[2mm]
0, \quad {\rm if \,\, not}. \\
\end{array}
\right.
\end{equation}
Then $\sqrt{|A_{\bf n}|} \leqslant {C}/{\bf n}$.
\end{proposition}

\noindent{\bf Proof of Proposition~\ref{proppi-1}.} Let $f^\ve_\mx = f^\ve_\mx(x, y, t)$
be a function that satisfies (\ref{rev-fff}).
Then we can write
\begin{equation}
\label{rev-fff-1}
\Vert f^\ve_\mx \Vert^2_{L^2(\Omega_T; L^2_{per}(Y_\mx))} =
\sum_{k=1}^{N_\ve} \big|\ve Y_\mx\big|\, \Vert \widehat f^\ve_{\mx,k}
\Vert^2_{L^2(0,T; L^2_{per}(Y_\mx))},
\end{equation}
where, due to (\ref{rev-fff}), we have that
\begin{equation}
\label{rev-fff-3}
\sum_{k=1}^{N_\ve} \big|\ve Y_\mx\big|\, \Vert \widehat f^\ve_{\mx,k}
\Vert^2_{L^2(0,T; L^2_{per}(Y_\mx))} \leqslant C^2.
\end{equation}
Now, for any ${\bf n} \in \mathbb{N}$ and $\ve >0$, let us introduce the set of "bad points" $A^\ve_{\bf n}$ defined by:
\begin{equation}
\label{rev-fff-4}
A^\ve_{\bf n} \eqdef \left\{x \in \Omega\,:\, \Vert \widehat f^\ve_{\mx,k(x,\ve)}
\Vert_{L^2(0,T; L^2_{per}(Y_\mx))} > {\bf n} \right\}.
\end{equation}
Let us estimate the measure of the set $A^\ve_{\bf n}$.
It follows from (\ref{rev-fff-3}) and (\ref{rev-fff-4}) that
$$
C^2 \geqslant \sum_{k=1}^{N_\ve} \big|\ve Y_\mx\big|\, \Vert \widehat f^\ve_{\mx,k(x,\ve)}
\Vert^2_{L^2(0,T; L^2_{per}(Y_\mx))} \geqslant
\sum_{x \in A^\ve_{\bf n}} \big|\ve Y_\mx\big|\, {\bf n}^2 = {\bf n}^2 \, | A^\ve_{\bf n}|.
$$
Therefore, $| A^\ve_{\bf n}|\leqslant C^2/{\bf n}^2$. By definition of limit inferior, for any
$\eta >0$ we have $A_{\bf n}\subseteq \liminf_{\ve\to 0} A^\ve_{{\bf n}-\eta}$,
(where $\ve$ denotes a {\it sequence} of real numbers).
Due to the continuity of the measure we get $|A_{\bf n}|\leqslant \liminf_{\ve\to 0} |A^\ve_{{\bf n}-\eta}|\leqslant C^2/({\bf n}-\eta)^2$.
Proposition \ref{proppi-1} is proved. \fin

We note that previously defined restricted functions are linked to ones appearing in
Proposition~\ref{proppi-1} by the following relation:
$$
f^\ve_{\mx,x_0}(y,t) = \widehat f^\ve_{\mx,k(x_0,\ve)}(y,t).
$$

In a similar way, taking into account the uniform estimate (\ref{dilo-12}),
we prove the following proposition.

\begin{proposition}
\label{proppi-2}
Let $f^\ve_\mx = f^\ve_\mx(x, y, t)$ be a dilated function such that
$\Vert f^\ve_\mx \Vert_{L^2(\Omega_T; H^{-1}_{per}(Y_\mx))} \leqslant C$
and let $B_{\bf n}$ be a set of points defined by
$$
B_{\bf n} \eqdef \left\{x \in \Omega\,:\, \liminf_{\ve\to0}\Vert \widehat f^\ve_{\mx,k(x,\ve)}
\Vert_{L^2(0,T; H^{-1}_{per}(Y_\mx))} \geqslant {\bf n} \right\},
$$
where the function $\widehat f^\ve_{\mx,k}$ is defined in (\ref{rev-fff-2}). Then
$\sqrt{|B_{\bf n}|} \leqslant C/{\bf n}$.
\end{proposition}

Now let us introduce ${\EuScript A}_{\bf n}$, the set of "bad points" for the
functions appearing in (\ref{dilo-8})-(\ref{dilo-10<2p}). We set:
$$
{\EuScript A}_{1,{\bf n}} \eqdef \left\{x \in \Omega\,:\,
\liminf_{\ve\to0}  \ve^{\theta/2 -1}\, \Vert \nabla_y \beta_\mx(s^\ve_{\mx,x})
\Vert_{L^2(0,T; L^2_{per}(Y_\mx))} \geqslant {\bf n} \right\};
$$
$$
{\EuScript A}_{2,{\bf n}} \eqdef \left\{x \in \Omega\,:\,
\liminf_{\ve\to0}\Vert p^\ve_{\mx,x}
\Vert_{L^2(0,T; L^2_{per}(Y_\mx))} \geqslant {\bf n} \right\};
$$
$$
{\EuScript A}_{3,{\bf n}} \eqdef \left\{x \in \Omega\,:\,
\liminf_{\ve\to0} \ve^{\theta/2 -1}\, \Vert \nabla_y p^\ve_{\mx,x}
\Vert_{L^2(0,T; L^2_{per}(Y_\mx))} \geqslant {\bf n} \right\};
$$
$$
{\EuScript A}_{4,{\bf n}} \eqdef \left\{x \in \Omega\,:\,
\liminf_{\ve\to0}\Vert \partial_t(\Phi_\mx\, s^\ve_{\mx,x})
\Vert_{L^2(0, T; H^{-1}_{per}(Y_\mx))} \geqslant {\bf n} \right\}.
$$
Here $s^\ve_{\mx,x}, p^\ve_{\mx,x}$ are defined in (\ref{new-smx}).
Then
\begin{equation}
\label{A_n}
{\EuScript A}_{\bf n} \eqdef \bigcup_{\ell=1}^4 {\EuScript A}_{\ell,{\bf n}}
\end{equation}
and, due to Propositions \ref{proppi-1}, \ref{proppi-2}, the measure of this set
satisfies the estimate $\sqrt{|{\EuScript A}_{\bf n}|} \leqslant {C}/{\bf n}$.

The following result holds.

\begin{lemma}
\label{lem-hitro-vyeb}
Let $s^\ve_{\mx,x_0}, p^\ve_{\mx,x_0}$ be the functions defined in (\ref{new-smx})
and $0 <\theta \leqslant 2$. Then for any $x_0 \in \Omega \setminus
{\EuScript A}_{\bf n}$, there is a subsequence $\ve = \ve_k$ still denoted by $\ve$ such that:
\begin{equation}
\label{dilo-8-vye}
0 \leqslant s^\ve_{\mx,x_0} \leqslant 1 \quad {\rm a.e.\,\, in\,\,} Y_\mx \times (0,T);
\end{equation}
\begin{equation}
\label{dilo-9-vye}
\Vert \nabla_y \beta_\mx(s^\ve_{\mx,x_0}) \Vert_{L^2(0, T; L^2_{per}(Y_\mx))} \leqslant C \ve^{1-\theta/2};
\end{equation}
\begin{equation}
\label{dilo-10-vye}
\Vert p^\ve_{\mx,x_0} \Vert_{L^2(0, T; L^2_{per}(Y_\mx))} \leqslant C;\quad
\Vert \nabla p^\ve_{\mx,x_0} \Vert_{L^2(0, T; L^2_{per}(Y_\mx))} \leqslant C\ve^{1-\theta/2};
\end{equation}
\begin{equation}
\label{dilo-12-vye}
\Vert \partial_t(\Phi_\mx\, s^\ve_{\mx,x_0}) \Vert_{L^2(0, T; H^{-1}_{per}(Y_\mx))}
\leqslant C,
\end{equation}
where $C = C({\bf n})$ is constant that does not depend on $x_0$ and $\ve$, and
 ${\bf n}$ is an arbitrary natural number.
\end{lemma}

\noindent{\bf Proof of Lemma \ref{lem-hitro-vyeb}.} First, we notice that the estimate
(\ref{dilo-8-vye}) follows immediately from (\ref{dilo-8}).
Let us prove, for example, (\ref{dilo-9-vye}). Taking into account that
$x_0 \in \Omega\setminus {\EuScript A}_{\bf n}$,
from the definition of the set ${\EuScript A}_{1,{\bf n}}$, we obtain
immediately the existence of a subsequence on which  (\ref{dilo-9-vye}) holds with constant $C$ depending only on ${\bf n}$.
The estimates (\ref{dilo-10-vye})-(\ref{dilo-12-vye}) are obtained
in a similar way. Lemma \ref{lem-hitro-vyeb} is proved. \fin

Using these estimates and applying Lemma 4.2 from \cite{our-siam}, we obtain
the following compactness result.

\begin{proposition}
\label{prop-vyeb}
Assume $0 < \theta \leqslant 2$. For any $x_0 \in \Omega\setminus {\EuScript A}_{\bf n}$,
on a subsequence extracted in Lemma~\ref{lem-hitro-vyeb}, the family $\{s^\ve_{\mx,x_0}\}_{\ve>0}$ is a
compact set in the space $L^q(Y_\mx \times (0,T))$ for all $q \in [1, \infty)$.
In the case $\theta < 2$ every  limit point of the sequence $\{s^\ve_{\mx,x_0}\}_{\ve>0}$
is independent of the fast variable $y$.
\end{proposition}

\section{Homogenization results}
\label{main-res}

\setcounter{equation}{0}

In this section we formulate and prove the main results of the paper corresponding to the
homogenized models for various rates of contrast. First, we introduce the
notation.

\begin{itemize}

\item[-] $S$, $P_w$, $P_n$ denote the homogenized wetting liquid saturation, wetting liquid pressure,
and nonwetting liquid pressure, respectively.

\item[-] $\Phi^\star = \Phi^\star(x)$ denotes the effective porosity and is given by:
\begin{equation}
\label{H-1}
\Phi^\star(x) \eqdef \Phi^{\rm H}_{\fr}(x)\, \frac{|Y_\fr|}{|Y_\mx|},
\end{equation}
where $\Phi^{\rm H}_{\fr}$ is defined in condition {(A.1)} and
$|Y_\ell|$ is the measure of the set $Y_\ell$ ($\ell = \fr, \mx$).

\item[-] $F^\star_w, F^\star_n$ denote the effective source terms and are given by:
\begin{equation}
\label{H-1FFF}
F^\star_w(x, t) \eqdef F^{\rm H}_{w}(x, t)\, \frac{|Y_\fr|}{|Y_\mx|} \quad {\rm and} \quad
F^\star_n(x, t) \eqdef F^{\rm H}_{n}(x, t)\, \frac{|Y_\fr|}{|Y_\mx|},
\end{equation}
where
\begin{equation}
\label{H-1FFF-source}
F^{\rm H}_w(x, t) \eqdef S^I_{\fr, w}\, f_I(x,t) - S\, f_P(x,t) \quad {\rm and}
\quad F^{\rm H}_n(x, t) \eqdef S^I_{\fr, n}\, f_I(x,t) - (1 - S)\, f_P(x,t)
\end{equation}
and where the functions $S^I_{\fr, w}, S^I_{\fr, n}, f_I, f_P$
are defined in (\ref{sour1}), (\ref{sour2}), respectively (see also (A.9)).

\item[-] $\mathbb{K}^\star = \mathbb{K}^\star(x)$ is the homogenized tensor with the entries
$\mathbb{K}^\star_{ij}$ defined by:
\begin{equation}
\label{H-2}
\mathbb{K}^{\star}_{ij}(x) \eqdef \frac{1}{|Y_\mx|}\, \int\limits_{Y_\fr}\, K(x, y)\,
\left[\nabla_y \xi_i + \vec e_i \right]\cdot\left[\nabla_y \xi_j + \vec e_j \right]\, dy,
\end{equation}
where $\xi_j = \xi_j(x, y)$ ($j = 1,\ldots,d$) is a $Y$-periodic solution to the auxiliary cell problem:
\begin{equation}
\label{H-20}
\left\{
\begin{array}[c]{ll}
- {\rm div}_y\, \big\{K(x, y) \nabla_y \xi_j\big\} = 0 \quad {\rm in} \,\, Y_{\fr}; \\[2mm]
\nabla_y \xi_j \cdot \vec \nu_y = - \vec e_j \cdot \vec \nu_y
\quad {\rm on} \,\, \Gamma_{\fr\mx};\\[2mm]
y \mapsto \xi_j(y)\quad Y-{\rm periodic}. \\
\end{array}
\right.
\end{equation}
\end{itemize}

\subsection{High contrast media: critical case, ${\boldsymbol \theta}{\bf=2}$}
\label{hig-crtic-subsec}

We study the asymptotic behavior of the solution to problem
(\ref{debut2}), (\ref{bc3})-(\ref{init1}) in the case $\varkappa(\ve) = \ve^2$ as $\ve \to 0$.
In particular, we are going to show that the effective model,
expressed in terms of the homogenized phase pressures,
reads:
\begin{equation}
\label{H-0}
\left\{
\begin{array}[c]{ll}
0 \leqslant S \leqslant 1 \quad {\rm in} \,\, \Omega_T; \\[2mm]
\displaystyle
\Phi^\star(x)\, \frac{\partial S}{\partial t}
- {\rm div}_x\, \bigg\{\mathbb{K}^\star(x)\, \lambda_{\,\fr,w}(S) \big(\nabla P_w - \vec g \big)  \bigg\}
=
{\EuScript Q}_w + F^\star_w \quad {\rm in} \,\, \Omega_T;
\\[5mm]
\displaystyle
- \Phi^\star(x)\, \frac{\partial S}{\partial t}
- {\rm div}_x\, \bigg\{\mathbb{K}^\star(x)\, \lambda_{\,\fr,n}(S)
\big(\nabla P_n - \vec g \big)  \bigg\}
=
{\EuScript Q}_n + F^\star_n \quad {\rm in} \,\, \Omega_T;\\[5mm]
P_{\fr,c}(S) = P_n - P_w \quad {\rm in} \,\, \Omega_T.
\end{array}
\right.
\end{equation}

For almost every point $x \in \Omega$ a matrix block $Y_\mx \subset \mathbb{R}^d$ is
suspended topologically. The system for flow in a matrix block is given by the so-called imbibition equation:
\begin{equation}
\label{H-4}
\left\{
\begin{array}[c]{ll}
\displaystyle
\Phi_m(y) \frac{\partial s}{\partial t}\, - {\rm div}_y \big\{ K(x, y) \nabla_y \beta_m(s)\big\} = 0
\quad {\rm in}\,\, Y_\mx \times \Omega_T; \\[3mm]
s(x, y, t) = {\EuScript P}(S(x,t)) \quad {\rm on}\,\, \Gamma_{\fr\mx} \times\Omega_T; \\[2mm]
s(x, y, 0) = S_m^{\,\bf 0}(x) \quad {\rm in}\,\, Y_\mx \times \Omega.\\
\end{array}
\right.
\end{equation}
Here $s$ denotes the wetting liquid saturation in the block $Y_\mx$ and
the function ${\EuScript P}(S)$ is defined by
\begin{equation}
\label{p-callig-moder}
{\EuScript P}(S)\, \eqdef\, (P_{c,m}^{-1} \circ P_{c,f})(S).
\end{equation}

For any $x \in \Omega$ and $t > 0$, the matrix-fracture sources are given by:
\begin{equation}
\label{H-6}
{\EuScript Q}_w \eqdef - \frac{1}{|Y_\mx|}\, \int\limits_{Y_\mx} \Phi_\mx(y)
\frac{\partial s}{\partial t}(x, y, t) \,dy = - {\EuScript Q}_n.
\end{equation}
The boundary conditions for the effective system (\ref{H-0}) are given by:
\begin{equation}
\label{H-7}
\left\{
\begin{array}[c]{ll}
P_{w} = P_{n} = 0 \quad {\rm on} \,\, \Gamma_{1} \times (0,T); \\[3mm]
\mathbb{K}^\star\,\lambda_{n}(S) \left(\nabla P_w - \vec g \right) \cdot \vec \nu =
\mathbb{K}^\star\, \lambda_{w}(S) (\nabla P_n - \vec g)\cdot \vec \nu = 0 \quad {\rm on} \,\,
\Gamma_{2} \times (0,T).\\
\end{array}
\right.
\end{equation}
Finally, the initial conditions read:
\begin{equation}
\label{H-9}
P_{w}(x, 0) = {\mathsf p}_{w}^{\bf 0}(x) \quad {\rm and}
\quad
P_{n}(x, 0) = {\mathsf p}_{n}^{\bf 0}(x)
\quad {\rm in} \,\, \Omega.
\end{equation}

The first main result of the paper is given by the following theorem.

\begin{theorem}
\label{t-hom-main}
Let $\varkappa(\ve) = \ve^2$ and let assumptions {(A.1)-(A.9)} be fulfilled.
Then the solution of the {\rm initial problem (\ref{debut2}), (\ref{inter-condit})-(\ref{init1})} converges
(up to a subsequence) in the two-scale sense to a weak solution of
the {\rm homogenized problem (\ref{H-0}), (\ref{H-4}), (\ref{H-6})-(\ref{H-9})}.
\end{theorem}

\noindent {\bf Proof of Theorem \ref{t-hom-main}.}
It is done in several steps. We start our analysis
by considering the system (\ref{debut2}). The main difficulty with the initial unknown functions
${\mathsf p}^\ve_{w}, {\mathsf p}^\ve_{n}$ in this system is that they do not possess the uniform
$H^1$-estimates (see Lemma \ref{lem-uniform}). It is important to notice that in the case of two-phase
incompressible flow it is possible to find appropriate but rather strong conditions which allow us to deal directly with the phase pressures in a space wider than $H^1$ (see \cite{yeh2}).
To overcome the difficulties appearing due to the absence of the uniform $H^1$-estimates, the
authors usually pass to the equivalent formulation of the problem in terms of the global pressure and
saturation. In our case it is done in subsection \ref{gp-relat} and the corresponding
weak formulation of the problem is then given in subsection \ref{def-weak-sol}.
Then using the convergence and compactness results from subsection \ref{ex-comp-sf} we pass to the limit
in equations (\ref{wf-1-gl-GP}), (\ref{gf-2-gl-GP}). This is done in subsections \ref{passage1} and
\ref{passage2}. In order, to pass to the {homogenized phase pressures} we make use
of the change of the unknown functions. Namely, we set, by the definition of the global pressure: $P_w \eqdef  {\mathsf P} + {\mathsf G}_{\fr,w}(S)$ and
$P_n \eqdef {\mathsf P} + {\mathsf G}_{\fr,n}(S)$. Then we rewrite the limit system
obtained in terms of the global pressure and saturation in terms of
the homogenized phase pressures. The passage to the limit in the matrix blocks makes use
of the dilation operator (see Section \ref{dil-oper} above).
Then we pass to the equivalent problem for the imbibition equation and, finally, obtain
system (\ref{H-4}).

\subsubsection{Passage to the limit in equation (\ref{wf-1-gl-GP})}
\label{passage1}

We set:
\begin{equation}
\label{paslim-1}
\varphi_w\left(x, \frac{x}{\ve}, t \right) \eqdef
\varphi(x, t) + \ve\, \zeta\left(x, \frac{x}{\ve}, t \right)
=
\varphi(x, t) + \ve\, \zeta_1(x, t)\, \zeta_2\left( \frac{x}{\ve}\right)
\eqdef \varphi(x, t) + \ve\, \zeta^\ve(x, t),
\end{equation}
where $\varphi \in {\EuScript D}(\Omega_T), \zeta_1 \in {\EuScript D}(\Omega_T),
\zeta_2 \in C^\infty_{per}(Y)$, and plug the function
$\varphi_w$ in (\ref{wf-1-gl-GP}). This yields:
$$
-\int\limits_{\Omega_{T}} {\bf 1}^\ve_\fr(x)\, \Phi^\ve_\fr(x)\, \widetilde S^\ve_\fr
\left[ \frac{\partial \varphi}{\partial t} +
\ve \frac{\partial \zeta^\ve}{\partial t} \right]\,\, dx\, dt
+
$$
$$
+
\int\limits_{\Omega_{T}}  {\bf 1}^\ve_\fr(x)\, K\left(x, \frac{x}{\ve}\right)
\bigg\{\lambda_{\fr,w} (\widetilde S^\ve_\fr)
\left(\nabla \widetilde {\mathsf P}^\ve_\fr - \vec g\right)
+ \nabla \beta_\fr(\widetilde S^\ve_\fr) \bigg\} \cdot \left[\nabla \varphi +
\ve \nabla_x \zeta^\ve + \nabla_y \zeta^\ve \right]\,\, dx\, dt
-
$$
$$
-\int\limits_{\Omega^\ve_{\mx,T}} \Phi_\mx\left(\frac{x}{\ve}\right) S^\ve_\mx
\left[ \frac{\partial \varphi}{\partial t} +
\ve \frac{\partial \zeta^\ve}{\partial t} \right]\,\, dx\, dt
+
$$
$$
+
\ve^2\, \int\limits_{\Omega^\ve_{\mx,T}} K\left(x, \frac{x}{\ve}\right) \bigg\{\lambda_{\mx,w} (S^\ve_\mx)
\left(\nabla {\mathsf P}^\ve_\mx - \vec g\right)
+ \nabla \beta_\mx(S^\ve_\mx) \bigg\} \cdot \left[\nabla \varphi +
\ve \nabla_x \zeta^\ve + \nabla_y \zeta^\ve \right]\,\, dx\, dt =
$$
\begin{equation}
\label{cr-10+1}
= \int\limits_{\Omega^\ve_{\fr,T}} \big(S^I_{\fr, w} f_I(x,t) - S^\ve_{\fr} f_P(x,t)\big)\,
\left[\varphi + \ve\, \zeta^\ve \right]\, dx dt.
\end{equation}

Taking into account Lemma~\ref{lem-uniform} and  the convergence results of Lemma \ref{2scale}
and Lemma \ref{conv-lemma-dil}, we pass to the limit in (\ref{cr-10+1}) as $\ve \to 0$ and
obtain the following homogenized equation:
$$
-\, |Y_\fr|\, \int\limits_{\Omega_T} \Phi^{\rm H}_\fr(x) S(x, t) \frac{\partial \varphi}{\partial t} \,dx\,dt +
$$
$$
+
\int\limits_{\Omega_T \times Y_\fr} K(x, y) \bigg\{ \lambda_{\,\fr,w}(S) \left[\nabla {\mathsf P}
+ \nabla_y {\mathsf w}_p - \vec g\right] + \nabla \beta_\fr(S) + \nabla_y {\mathsf w}_s \bigg\} \cdot \left[\nabla \varphi + \zeta_1 \nabla_y \zeta_2 \right] \,dy\, dx\, dt
=
$$
\begin{equation}
\label{fin-cut-1}
=
\int\limits_{\Omega_T \times Y_\mx} \Phi_\mx(y)\, s(x, y, t) \, \frac{\partial \varphi}{\partial t}
\,\,dy\, dx\, dt + |Y_\fr|\, \int\limits_{\Omega_T} F^{\rm H}_w\, \varphi\, dx\, dt,
\end{equation}
where $F^{\rm H}_w$ is given by (\ref{H-1FFF-source}).

\subsubsection{Passage to the limit in equation (\ref{gf-2-gl-GP})}
\label{passage2}

Equation  (\ref{gf-2-gl-GP}) is treated in the same way as equation  (\ref{wf-1-gl-GP}). Taking the test function of
the form (\ref{paslim-1}) and using the same arguments  we can pass to a limit $\ve \to 0$ and obtain  the following homogenized equation:
$$
|Y_\fr|\, \int\limits_{\Omega_T} \Phi^{\rm H}_\fr(x) S(x, t)
\frac{\partial \varphi}{\partial t} \,dx\,dt +
$$
$$
+
\int\limits_{\Omega_T \times Y_\fr} K(x, y) \bigg\{ \lambda_{\fr,n}(S) \left[\nabla {\mathsf P}
+ \nabla_y {\mathsf w}_p - \vec g\right] - \nabla \beta_\fr(S) - \nabla_y {\mathsf w}_s \bigg\}
\cdot \left[\nabla \varphi + \zeta_1 \nabla_y \zeta_2 \right] \,dy\, dx\, dt
=
$$
\begin{equation}
\label{fin-cut-2}
= -
\int\limits_{\Omega_T \times Y_\mx} \Phi_\mx(y)\, s(x, y, t) \, \frac{\partial \varphi}{\partial t}
\,\,dy\, dx\, dt + |Y_\fr|\, \int\limits_{\Omega_T} F^{\rm H}_n\, \varphi\, dx\, dt.
\end{equation}

\subsubsection{Identification of the corrector functions ${\mathsf w}_p$, ${\mathsf w}_s$ and homogenized equations}
\label{ident-p-beta}

In this section we identify the corrector functions ${\mathsf w}_p$, ${\mathsf w}_s$ appearing
in the equations (\ref{fin-cut-1}), (\ref{fin-cut-2}) and obtain the desired homogenized
system (\ref{H-0}).

Consider the equations (\ref{fin-cut-1}), (\ref{fin-cut-2}). Setting $\varphi \equiv 0$, we get:
\begin{equation}
\label{hernya-3}
\int\limits_{Y_{\fr}} K(x, y)\,\bigg\{ \lambda_{\,\fr,w}(S)
\big[\nabla P + \nabla_{y} {\mathsf w}_p - \vec{g} \big]
+ \big[\nabla\beta_\fr + \nabla_{y} {\mathsf w}_s\big] \bigg\} \cdot \nabla_y \zeta_2(y)\,dy = 0
\end{equation}
and
\begin{equation}
\label{hernya-4}
\int\limits_{Y_{\fr}} K(x, y)\,\bigg\{ \lambda_{\,\fr,n}(S)
\big[\nabla P + \nabla_{y} {\mathsf w}_p - \vec{g} \big]
- \big[\nabla\beta_\fr + \nabla_{y} {\mathsf w}_s\big] \bigg\} \cdot \nabla_y \zeta_2(y)\,dy = 0.
\end{equation}
Now adding (\ref{hernya-3}) and (\ref{hernya-4}) and taking into account condition
{(A.4)} and the fact that the saturation $S$ does not depend on the fast variable $y$, we obtain:
\begin{equation}
\label{hernya-5}
\int\limits_{Y_{\fr}} K(x, y)\,\bigg\{\nabla P + \nabla_{y} {\mathsf w}_p  -
\vec{g} \bigg\} \cdot \nabla_y \zeta_2(y)\,dy = 0.
\end{equation}
Then we proceed in a standard way (see, e.g., \cite{hor}). Let $\xi_j = \xi_j(x, y)$
($j = 1,..,d$) be the $Y$-periodic solution of the auxiliary cell problem (\ref{H-20}).
Then the function ${\mathsf w}_p$ can be represented as:
\begin{equation}
\label{1bvp-2}
{\mathsf w}_p(x, y, t) = \sum^d_{j=1} \xi_j(x, y) \left[\frac{\partial\, P}{\partial x_j}(x, t) - g_j\right].
\end{equation}

Now we turn to the identification of the function ${\mathsf w}_s$. From (\ref{hernya-3})
and (\ref{hernya-5}), we get:
\begin{equation}
\label{wbet-1}
\int\limits_{Y_{\fr}} K(x, y)\, \bigg\{ \nabla\beta_{\fr} + \nabla_{y} {\mathsf w}_s \bigg\} \cdot \nabla_y \zeta_2(y)\,dy = 0.
\end{equation}
Then as in the previous case, we obtain that
\begin{equation}
\label{wbet-2}
{\mathsf w}_s(x, y, t) = \sum^d_{j=1} \xi_j(x, y) \frac{\partial\, \beta_\fr(S)}{\partial x_j}(x, t).
\end{equation}

\subsubsection{Effective equations in terms of the global pressure and saturation}

We start by obtaining the corresponding homogenized equation for the wetting phase.
Choosing $\zeta_2 = 0$ in    (\ref{fin-cut-1}), we get:
$$
\Phi^\star(x) \, \dfrac{\partial S}{\partial t} -
{\rm div}_x \bigg\{\mathbb{K}^\star(x)\,
\big[\lambda_{\,\fr,w}(S)\,\nabla P + \nabla\beta_{\fr}(S) - \lambda_{\,\fr,w}(S)\, \vec{g}\big]
\bigg\} =
$$
\begin{equation}
=
- \frac{1}{|Y_{m}|}\, \int\limits_{Y_{\mx}} \Phi_{m}(y)
\dfrac{\partial s}{\partial t}(x, y, t)\,dy + F^\star_w(x, t),
\label{1newpoint-7:2ok}
\end{equation}
where the effective porosity $\Phi^\star$, the effective source term $F^\star_w$, and
the homogenized permeability tensor $\mathbb{K}^\star$
are defined in \eqref{H-1}, \eqref{H-1FFF} and
\eqref{H-2}, respectively.

In a similar way, choosing $\zeta_2 = 0$ in
 equation (\ref{fin-cut-2}), we derive the second homogenized equation:
$$
- \Phi^\star(x) \, \dfrac{\partial S}{\partial t} -
{\rm div}_x \bigg\{\mathbb{K}^\star(x)\,
\big[\lambda_{\,\fr,n}(S)\,\nabla P + \nabla\beta_{\fr}(S) - \lambda_{\,\fr,n}(S)\, \vec{g}\big]
\bigg\} =
$$
\begin{equation}
=
\frac{1}{|Y_{m}|}\, \int\limits_{Y_{\mx}} \Phi_{m}(y)
\dfrac{\partial s}{\partial t}(x, y, t)\,dy + F^\star_n(x, t),
\label{hernya-18ok}
\end{equation}
where $F^\star_n$ denotes the effective source term defined in \eqref{H-1FFF}.

\subsubsection{Effective equations in terms of the phase pressures}
\label{hom-eq-in-ph-pres}

Let us introduce now the functions that is naturally to call the {homogenized phase pressures}.
Namely, we set, by the definition:
\begin{equation}
\label{phase-pres}
P_{w} \eqdef {\mathsf P} + {\mathsf G}_{\fr,w}(S) \quad {\rm and}
\quad P_n \eqdef {\mathsf P} + {\mathsf G}_{\fr,n}(S),
\end{equation}
where the functions ${\mathsf G}_{\fr,w}, {\mathsf G}_{\fr,n}$ are defined in Section \ref{gp-relat}.
Then it easy to see that the homogenized equations can be rewritten as follows:
\begin{equation}
\label{final-h-press}
\left\{
\begin{array}[c]{ll}
\displaystyle
\Phi^\star(x)\, \frac{\partial S}{\partial t}
- {\rm div}_x\, \bigg\{\mathbb{K}^\star(x)\, \lambda_{\,\fr,w}(S) \big(\nabla P_w - \vec g \big)  \bigg\}
=
{\EuScript Q}_w + F^\star_w \quad {\rm in} \,\, \Omega_T;
\\[5mm]
\displaystyle
- \Phi^\star(x)\, \frac{\partial S}{\partial t}
- {\rm div}_x\, \bigg\{\mathbb{K}^\star(x)\, \lambda_{\,\fr,n}(S)
\big(\nabla P_n - \vec g \big)  \bigg\}
=
{\EuScript Q}_n + F^\star_n \quad {\rm in} \,\, \Omega_T;\\[5mm]
P_c(S) = P_n - P_w \quad {\rm in} \,\, \Omega_T.
\end{array}
\right.
\end{equation}

\subsubsection{Flow equations in the matrix block}
\label{identif-l-theta}

In this section, following the ideas of the papers \cite{ba-lp-doubpor,blm,choq}, we
obtain the system (\ref{H-4}) describing the behavior of the function $s$ which is
involved in the definition of the matrix-fracture source term.
Briefly, we pass to the limit in the equations for the dilated functions
for fixed $k$ and then by density arguments the limit equations will be
obtained. We recall that the equations for the dilated functions are
already obtained in Lemma \ref{lem-dil} from subsection \ref{dil-func-prop}. Namely,
for almost all $x \in \Omega$, the functions $s^\ve_\mx, p^\ve_\mx$
satisfy the following variational problem:

\noindent{for all $\phi_n,\phi_w \in L^2(0,T;H^1_0(Y_{\mx}))\cap H^1(0,T; L^2(Y_{\mx}))$, $\phi_n(T)=\phi_w(T) =0$,}
$$
- \int\limits_0^T\int\limits_{Y_{\mx}} \Phi_\mx(y) s^\ve_\mx  \frac{\partial \phi_w}{\partial t}\, dy
- \int\limits_0^T\int\limits_{Y_{\mx}} \Phi_\mx(y) (\mathfrak{D}^\ve S^{\bf 0}_\mx)
\frac{\partial \phi_w}{\partial t}(0)\, dy
$$
\begin{equation}
\label{ebio-ident-1v-w}
+  \int\limits_0^T\int\limits_{Y_{\mx}} \bigg\{K(x, y) \left[
\lambda_{\mx,w} (s^\ve_\mx) \nabla_y p^\ve_\mx +
\nabla_y \beta_\mx(s^\ve_\mx) - \ve\,\lambda_{\mx,w} (s^\ve_\mx) \vec{g}\, \right] \bigg\}
\cdot \nabla_y \phi_w \, dy = 0;
\end{equation}
$$
\int\limits_0^T\int\limits_{Y_{\mx}}  \Phi_\mx(y) s^\ve_\mx \frac{\partial \phi_n }{\partial t}\, dy
+ \int\limits_0^T\int\limits_{Y_{\mx}}  \Phi_\mx(y)  (\mathfrak{D}^\ve S^{\bf 0}_\mx)
\frac{\partial \phi_n }{\partial t}\, dy
$$
\begin{equation}
\label{ebio-ident-1v-n}
+\int\limits_0^T\int\limits_{Y_{\mx}} \bigg\{K(x, y) \left[ \lambda_{\mx,n} (s^\ve_\mx) \nabla_y p^\ve_\mx -
\nabla_y \beta_\mx(s^\ve_\mx) - \ve\, \lambda_{\mx,n}(s^\ve_\mx)\,\vec{g}\,
\right] \bigg\} \cdot \nabla_y \phi_n \, dy= 0
\end{equation}
{with the boundary conditions (\ref{dilo-6})}.

The uniform estimates for the functions $s^\ve_\mx, p^\ve_\mx$ imply the convergence results
of $\langle s^\ve_\mx, p^\ve_\mx \rangle$ to  $\langle s, p \rangle$
in a weak sense (see Lemma \ref{conv-lemma-dil}). Thus, the limit behavior of the dilated
functions $s^\ve_\mx$, $p^\ve_\mx$ is determined.
However, the convergence results of
Lemma \ref{conv-lemma-dil} are not sufficient for derivation of the equations for the
limit functions $\langle s, p \rangle$. To overcome this difficulty, in
Section \ref{dil-func-conv} we pass to the restrictions of the functions
$s^\ve_\mx, p^\ve_\mx$ to $K^\ve_{x_0}$ defined in (\ref{new-smx}).
Evidently, they are constants in the slow variable $x$.
Introducing the set ${\EuScript A}_{\bf n}$ of "bad points" (\ref{A_n}), by
Lemma~\ref{lem-hitro-vyeb} we have the uniform estimates (\ref{dilo-8-vye})-(\ref{dilo-12-vye})
for the functions $s^\ve_{\mx,x_0}, p^\ve_{\mx,x_0}$. For any $\ve > 0$, the pair of functions
$\langle s^\ve_{\mx,x_0}, p^\ve_{\mx,x_0} \rangle$ is a solution to
problem (\ref{ebio-ident-1v-w}), (\ref{ebio-ident-1v-n}) in $Y_\mx \times (0,T)$. Moreover, the
compactness result, i.e., Proposition
\ref{prop-vyeb} is established for the family $\{s^\ve_{\mx,x_0}\}_{\ve>0}$.
Having established these results, we are in position to complete the proof
of Theorem \ref{t-hom-main}. The uniform estimates for the functions
$s^\ve_{\mx,x_0}, p^\ve_{\mx,x_0}$ from Lemma \ref{lem-hitro-vyeb} and
the compactness result formulated in Proposition \ref{prop-vyeb}
allow us to obtain the following convergence results.

\begin{lemma}
\label{conv-k-dil}
Let $x_0\in \Omega\setminus {\EuScript A}_{\bf n}$.
There exist functions $s_{x_0}, p_{x_0}$, and $\beta_\mx(s_{x_0})$ such that
up to a subsequence:
\begin{equation}
\label{2s-k-1}
s^\ve_{\mx,x_0} \to s_{x_0} \,\, {\rm strongly\,\, in}\,\, L^q(Y_\mx \times (0,T))\,\,
\forall \ 1 \leqslant q < +\infty;
\end{equation}
\begin{equation}
\label{2s-k-2}
p^\ve_{\mx,x_0} \rightharpoonup p_{x_0} \,\, {\rm weakly\,\, in}\,\, L^2(0, T; H^1_{per}(Y_\mx));
\end{equation}
\begin{equation}
\label{2s-k-40}
\beta_\mx(s^\ve_{\mx,x_0}) \rightharpoonup \beta_\mx(s_{x_0}) \,\, {\rm weakly\,\, in}\,\,
L^2(0, T; H^1_{per}(Y_\mx));
\end{equation}
\begin{equation}
\label{2s-k-4}
\beta_\mx(s^\ve_{\mx,x_0}) \to \beta_\mx(s_{x_0}) \,\, {\rm strongly\,\, in}\,\,
L^q(Y_\mx \times (0,T))\,\, \forall \ 1 \leqslant q < +\infty;
\end{equation}
\begin{equation}
\label{2s-k-4-front}
\beta_\mx(s^\ve_{\mx,x_0})\big|_{\Gamma_{\mx\fr}} \to \beta_\mx(s_{x_0})\big|_{\Gamma_{\mx\fr}}
\,\, {\rm weakly\,\, in}\,\, L^2(0, T; L^2(\Gamma_{\mx\fr}));
\end{equation}
\begin{equation}
\label{2s-k-5-front}
p^\ve_{\mx,x_0}\big|_{\Gamma_{\mx\fr}} \to p_{x_0}\big|_{\Gamma_{\mx\fr}}  \,\,
{\rm weakly\,\, in}\,\, L^2(0, T; L^2(\Gamma_{\mx\fr})).
\end{equation}

\end{lemma}

As in subsections \ref{passage1}, \ref{passage2} we can easily pass to the limit in
(\ref{ebio-ident-1v-w}) and (\ref{ebio-ident-1v-n}). We get the following system of equations:
\begin{align}
\label{ebio-ident-3v-w}
- \int\limits_0^T\int\limits_{Y_{\mx}} \Phi_\mx(y) &s_{x_0}  \frac{\partial \phi_w}{\partial t}\, dy
- \int\limits_0^T\int\limits_{Y_{\mx}} \Phi_\mx(y) S^{\bf 0}_{x_0}
\frac{\partial \phi_w}{\partial t}(0)\, dy \\
&+  \int\limits_0^T\int\limits_{Y_{\mx}} \bigg\{K(x, y) \big[
\lambda_{\mx,w} (s_{x_0}) \nabla_y p_{x_0} +
\nabla_y \beta_\mx(s_{x_0})\big] \bigg\}\cdot \nabla_y \phi_w \, dy
= 0; \nonumber \\
\int\limits_0^T\int\limits_{Y_{\mx}}  \Phi_\mx(y) & s_{x_0} \frac{\partial \phi_n }{\partial t}\, dy
+ \int\limits_0^T\int\limits_{Y_{\mx}}  \Phi_\mx(y)  S^{\bf 0}_{x_0}
\frac{\partial \phi_n }{\partial t}\, dy
\label{ebio-ident-3v-n} \\
&+\int\limits_0^T\int\limits_{Y_{\mx}} \bigg\{K(x, y)
\big[ \lambda_{\mx,n} (s_{x_0}) \nabla_y p_{x_0} -
\nabla_y \beta_\mx(s_{x_0}) \big] \bigg\} \cdot \nabla_y \phi_n \, dy= 0,\nonumber
\end{align}
where we have used the fact that $\mathfrak{D}^\ve S^{\bf 0}_{x_0}\to S^{\bf 0}_{x_0}$ strongly in $L^2(Y_m)$
for almost all $x_0\in \Omega$.

Now we turn to the boundary condition for $s_{x_0}$ on $\Gamma_{\mx\fr}$.
From Corollary~\ref{corol-dilop} we know that for a.e. $x_0$,
\begin{align*}
{\cal M}(\beta_\fr(\mathfrak{D}^\ve(\widetilde S^\ve_\fr(x_0,\cdot,\cdot)))) \to
{\cal M}(\beta_\fr(S(x_0,\cdot,\cdot)))\quad\text{strongly in }\;  L^2(0,T; L^2(\Gamma_{\mx\fr})),
\end{align*}
where ${\cal M}$ is the function given in (\ref{function:M}). Therefore, for a.e.
$x_0 \in \Omega\setminus {\EuScript A}_{\bf n}$, from (\ref{dilo-6}) and (\ref{2s-k-4-front})
we have:
\begin{align*}
\beta_\mx(s_{x_0})\big|_{\Gamma_{\mx\fr}} =
{\cal M}(\beta_\fr(S(x_0,\cdot,\cdot)))\big|_{\Gamma_{\mx\fr}},
\end{align*}
or, equivalently
\begin{align}
s_{x_0} = {\EuScript P}(S(x_0,\cdot)) \quad \text{on }\; \Gamma_{\mx\fr}\times (0,T).
\label{BC:matrix:eff}
\end{align}
Note also that  it follows  from (\ref{dilo-6}) that the convergence in (\ref{2s-k-4-front})
is strong in $L^2(0, T; L^2(\Gamma_{\mx\fr}))$. This, together with convergence
(\ref{2s-k-5-front}) and Lipschitz continuity  of the functions ${\mathsf G}_{\ell,g} $,
${\mathsf G}_{\ell,w} $, enables us to pass to the limit in the boundary condition
for dilated global pressure (\ref{dilo-6a}) using the two-scale convergence on
$\Gamma_{\fr\mx}$, and get
\begin{align}
p_{x_0} + {\mathsf G}_{\mx, w}(s_{x_0})={\mathsf P}(x_0, \cdot)+ {\mathsf G}_{\fr, w}(S(x_0,\cdot))
\quad \text{on }\; \Gamma_{\mx\fr}\times (0,T).
\label{BC:matrix:ess:P}
\end{align}
In the same way we also get
\begin{align}
p_{x_0} + {\mathsf G}_{\mx, g}(s_{x_0})={\mathsf P}(x_0, \cdot)+ {\mathsf G}_{\fr, g}(S(x_0,\cdot))
\quad \text{on }\; \Gamma_{\mx\fr}\times (0,T).
\label{BC:matrix:ess:P-g}
\end{align}

Thus the system which is satisfied by the limit $\langle s_{x_0}, p_{x_0}\rangle$ is obtained for
any $x_0 \in \Omega\setminus {\EuScript A}_{\bf n}$. Now it remains to make the link
between the functions $s_{x_0}, p_{x_0}$ and the limits $s, p$ of the sequences
$\{s^\ve_\mx\}_{\ve>0}$, $\{p^\ve_\mx\}_{\ve>0}$. First, we observe that the
convergent subsequence in Lemma~\ref{conv-k-dil} depends on  point
$x_0 \not\in {\EuScript A}_{\bf n}$. To avoid this difficulty we will prove
(see subsection \ref{identif-l-theta} below) that the problem  (\ref{ebio-ident-3v-w}),
(\ref{ebio-ident-3v-n}) with the corresponding boundary conditions (\ref{BC:matrix:eff}),
(\ref{BC:matrix:ess:P}), and (\ref{BC:matrix:ess:P-g}) has a unique weak solution.
Then the convergence results from Lemma~\ref{conv-k-dil} hold for the whole
sequences, as $\ve\to0$. Since the functions
$s_{x_0} = s(x_{x_0}, y, t), p_{x_0} = p(x_{x_0}, y, t)$ satisfy  (\ref{ebio-ident-3v-w})-(\ref{BC:matrix:ess:P-g})
for almost all $x_0 \in \Omega \setminus {\EuScript A}_{\bf n}$, we conclude that
$s$ and  $p$ are weak solution of the following system of equations:
\begin{equation}
\label{ebio-ident-4}
\left\{
\begin{array}[c]{ll}
0 \leqslant s \leqslant 1 \quad {\rm in} \,\, Y_\mx\times\Omega_T; \\[2mm]
\displaystyle
\Phi_\mx(y) \frac{\partial s}{\partial t} -
{\rm div}_y\, \bigg\{K(x, y) \left[
\lambda_{\mx,w} (s) \nabla_y p + \nabla_y \beta_\mx(s) \right] \bigg\}
= 0 \quad {\rm in} \,\, Y_\mx\times\Omega_T; \\[5mm]
\displaystyle
- \Phi_\mx(y) \frac{\partial s}{\partial t}
-
{\rm div}_y\, \bigg\{K(x, y) \left[ \lambda_{\mx,n} (s)
\nabla_y p - \nabla_y \beta_\mx(s) \right] \bigg\} = 0
\quad {\rm in} \,\, Y_\mx\times\Omega_T.\\[2mm]
\end{array}
\right.
\end{equation}

The system is completed by the corresponding boundary and initial conditions:
\begin{equation}
\label{last-bc-s}
\left\{
\begin{array}[c]{ll}
{\mathsf P} + {\mathsf G}_{\fr,w}(S) =
p + {\mathsf G}_{\mx,w}(s) \quad {\rm on}\,\, \Gamma_{\fr\mx} \times\Omega_T;\\[4mm]
{\mathsf P} + {\mathsf G}_{\fr,n}(S) =
p + {\mathsf G}_{\mx,n}(s) \quad {\rm on}\,\, \Gamma_{\fr\mx} \times\Omega_T; \\[4mm]
s(x, y, t) = {\EuScript P}(S(x,t)) \quad {\rm on}\,\, \Gamma_{\fr\mx} \times\Omega_T ,\\[2mm]
s(x, y, 0) = S^{\,\bf 0}(x) \quad {\rm in}\,\,  Y_\mx \times \Omega.
\end{array}
\right.
\end{equation}

Thus, we have identified $s$ and $p$ for $x\in  \Omega \setminus {\EuScript A}_{\bf n}$. Since
by Propositions \ref{proppi-1}, \ref{proppi-2}, the measure of the set ${\EuScript A}_{\bf n}$
goes to zero as ${\bf n}\to \infty$ we conclude that our conclusion holds a.e. in $\Omega$.

The proof of the uniqueness of the solution to problem (\ref{ebio-ident-4}) will be done as follows. First, we reduce
the system (\ref{ebio-ident-4}) to a boundary value problem for the so-called
{imbibition equation} and then make use of the uniqueness result from \cite{vazquez}.
Equation (\ref{H-4})$_1$ is the well known {generalized porous medium equation}
(see, e.g., \cite{vazquez}).

\begin{lemma}
\label{lemma-der-imb-eqn}
Let $s = s(x, y, t)$ be the solution of the cell problem
\eqref{ebio-ident-4}-\eqref{last-bc-s}. Then $s$ satisfies the boundary value problem (\ref{H-4}).
\end{lemma}

\noindent{\bf Proof of Lemma \ref{lemma-der-imb-eqn}.} First we observe that
it follows from the boundary conditions (\ref{last-bc-s}) that the function $s$ does not depend on
$y$ on $\Gamma_{\fr\mx} \times\Omega_T$. Then the global pressure $p$ does not depend on $y$ on
$\Gamma_{\fr\mx} \times\Omega_T$. Namely, we can write that
\begin{equation}
\label{H-9-then-2}
p(x, y, t) = p_{\Gamma}(x, t) \quad {\rm on} \,\, \Gamma_{\fr\mx} \times\Omega_T.
\end{equation}
By summing the two equations in \eqref{ebio-ident-4} we get:
\begin{equation}
\label{H-18+}
- {\rm div}\, \big\{K(x, y)\, \lambda_{\mx}(s) \nabla_y p \big\} = 0 \quad {\rm in}\,\,
Y_\mx \times \Omega.
\end{equation}
Then multiplying the equation \eqref{H-18+} by $(p -  p_{\Gamma})$
and integrating over $Y_\mx \times \Omega_T$, using \eqref{H-9-then-2} and conditions {(A.2)},
{(A.4)} we obtain:
$$
0 = \int\limits_{Y_\mx \times \Omega_T} K(x, y)\,
\lambda_{\mx}(s) \nabla_y p \cdot \nabla_y p \,\, dx\,dy\,dt \geqslant
k_{\rm min}\,\, L_0\, \int\limits_{Y_\mx \times \Omega_T} |\nabla_y p|^2 \,\, dx\,dy\,dt,
$$
which gives
$\nabla_y p = 0 \quad {\rm a.e.\,\, in} \,\, Y_\mx \times \Omega_T$.
This result allows us to reduce the two equations in the problem \eqref{ebio-ident-4} to only one,
as announced in (\ref{H-4}). This completes the proof of Lemma~\ref{lemma-der-imb-eqn}.
\fin

Now we turn to the proof of the uniqueness of the solution to (\ref{H-4}). This
proof is given in Theorem 5.3 from \cite{vazquez}. For reader's convenience we
discuss it briefly in the following lemma.

\begin{lemma}
\label{uniqueness-imbib}
Under our standing assumptions, there is a unique weak solution to problem (\ref{H-4}).
\end{lemma}

\noindent{\bf Proof of Lemma \ref{uniqueness-imbib}.} First, we introduce the weak
formulation of problem (\ref{H-4}). Omitting, for the sake of simplicity, the dependence
on the slow variable $x$, we have:

{\it for any function $\eta \in C^1(\overline{Y}_{\mx,T})$, where $Y_{\mx,T} \eqdef Y_\mx \times (0,T)$,
vanishing on $\Gamma_{\fr\mx}$ and such that $\eta(x, T) = 0$},
\begin{equation}
\label{uni-im-1}
\int\limits_{Y_{\mx,T}} \left\{K(y)\, \nabla_y \beta_\mx(s) \cdot \nabla_y \eta -
s\, \frac{\partial \eta}{\partial t} \right\}\, dy\, dt =
\int\limits_{Y_\mx}  S_m^{\bf 0}(x)\, \eta(y, 0)\, dy,
\end{equation}
.

Suppose now that we have two solutions $s_1$ and $s_2$ satisfying (\ref{uni-im-1}). Then
denoting $W_i \eqdef \beta_\mx(s_i)$, from (\ref{uni-im-1}), we have:
\begin{equation}
\label{uni-im-2}
\int\limits_{Y_{\mx,T}} \left\{K(y)\, \nabla_y (W_1 - W_2) \cdot \nabla_y \eta -
(s_1 - s_2)\, \frac{\partial \eta}{\partial t} \right\}\, dy\, dt = 0
\end{equation}
for all $\eta$. Then we use as a special test function $\eta = \widehat\eta$, see e.g. \cite{vazquez}:
\begin{equation}
\label{uni-im-3}
\widehat\eta \eqdef \left\{
\begin{array}[c]{ll}
\displaystyle
\int_t^T \big[W_1(x, \varsigma) - W_2(x, \varsigma)\big]\, d\varsigma \quad
{\rm if}\,\, 0 < t < T;\\[5mm]
0, \quad {\rm if}\,\, t \geqslant T.
\end{array}
\right.
\end{equation}
Then, plugging (\ref{uni-im-3}) in (\ref{uni-im-2}), we get:
\begin{equation}
\label{uni-im-4}
\int\limits_{Y_{\mx,T}} (s_1 - s_2)\, (W_1 - W_2) \, dy\, dt +
\int\limits_{Y_{\mx,T}} K(y)\, \nabla_y (W_1 - W_2) \cdot
\left\{ \int\limits_t^T \nabla_y (W_1 - W_2)\, d\varsigma \right\}\, dy\, dt  = 0.
\end{equation}
Integration of the last term leads to the following relation:
\begin{equation}
\label{uni-im-5}
\int\limits_{Y_{\mx,T}} (s_1 - s_2)\, \big(\beta_\mx(s_1) - \beta_\mx(s_2)\big) \, dy\, dt +
\frac12\, \int\limits_{Y_\mx} K(y)\,
\left[\int\limits_0^T \nabla_y \big(\beta_\mx(s_1) - \beta_\mx(s_2)\big)\, d\varsigma \right]^2\, dy = 0.
\end{equation}
Due to the monotonicity of the function $\beta_\mx$, the first term in (\ref{uni-im-2}) is
non-negative. Therefore we can conclude that $s_1 = s_2$ a.e. in $Y_{\mx,T}$.
Lemma \ref{uniqueness-imbib} is proved. \fin

This completes the proof of Theorem \ref{t-hom-main}. \fin

\subsection{Very high contrast media: ${\boldsymbol \theta}{\bf>2}$}
\label{very-hig-crtic-subsec}

We study the asymptotic behavior of the solution to problem
(\ref{debut2}), (\ref{bc3})-(\ref{init1}) as $\ve \to 0$ in the case $\varkappa(\ve) = \ve^\theta$ with $\theta > 2$.
In particular, we are going to show that the effective model reads:
\begin{equation}
\label{H-0-veryhigh}
\left\{
\begin{array}[c]{ll}
0 \leqslant S \leqslant 1 \quad {\rm in} \,\, \Omega_T; \\[2mm]
\displaystyle
\Phi^\star(x)\, \frac{\partial S}{\partial t}
- {\rm div}_x\, \bigg\{\mathbb{K}^\star(x)\, \lambda_{\,\fr,w}(S) \big(\nabla P_w - \vec g \big)  \bigg\}
=
F^\star_w \quad {\rm in} \,\, \Omega_T;
\\[5mm]
\displaystyle
- \Phi^\star(x)\, \frac{\partial S}{\partial t}
- {\rm div}_x\, \bigg\{\mathbb{K}^\star(x)\, \lambda_{\,\fr,n}(S)
\big(\nabla P_n - \vec g \big)  \bigg\}
=
F^\star_n \quad {\rm in} \,\, \Omega_T;\\[5mm]
P_{\fr,c}(S) = P_n - P_w \quad {\rm in} \,\, \Omega_T,
\end{array}
\right.
\end{equation}
where the effective porosity $\Phi^\star$, the effective source terms $F^\star_w, F^\star_n$,
and the homogenized permeability tensor $\mathbb{K}^\star$
in \eqref{H-0-veryhigh} are defined in \eqref{H-1}, \eqref{H-1FFF} and
\eqref{H-2}, respectively.

The boundary and the initial conditions for the system (\ref{H-0-veryhigh}) are given by (\ref{H-7}), (\ref{H-9}).

We see that in this case the matrix blocks have a vanishing, as $\ve \to 0$, influence
on the effective flow. This means that in the case of very high contrast, the medium
behaves as a perforated one.

The second main result of the paper is as follows.

\begin{theorem}
\label{t-hom-main>2}
Let $\varkappa(\ve) = \ve^\theta$ with $\theta > 2$ and let assumptions
{(A.1)-(A.9)} be fulfilled. Then the solution of the {\rm initial problem
(\ref{debut2}), (\ref{inter-condit})-(\ref{init1})} converges (up to a subsequence) in
the two-scale sense to a weak solution of the
{\rm homogenized problem \eqref{H-0-veryhigh}, (\ref{H-7}), (\ref{H-9})}.
\end{theorem}

\noindent{\bf Proof of Theorem \ref{t-hom-main>2}.}
Let $\theta>2$. In the proof of Theorem \ref{t-hom-main>2} we follow
the lines of the proof of Theorem \ref{t-hom-main}. Namely, arguing as
in Sections \ref{passage1}, \ref{passage2}, \ref{ident-p-beta}, we obtain
the homogenized equations (\ref{1newpoint-7:2ok}), (\ref{hernya-18ok}).

Now we want to show that in the case of the very high contrast, the
model behaves as in the perforated media, i.e., the matrix blocks
are totally impermeable and the additional matrix-source term equals
zero. As in the paper \cite{yeh2}, we prove the following result.

\begin{lemma}
\label{nosource}
The following equation holds true:
$$
\Phi_\mx(y)\, \frac{\partial s}{\partial t}(x, y, t) = 0 \quad
{\rm in} \,\,\, Y_\mx \times \Omega_T.
$$
\end{lemma}

\noindent{\bf Proof of Lemma \ref{nosource}.}
Let us define the function:
$$
{\EuScript F}^{\,\ve}(x,t) \eqdef \ve^{\frac{\theta}{2}}\,K^\ve(x)\, \bigg\{\lambda_{\mx,w} (S^\ve_\mx)
\left(\nabla {\mathsf P}^\ve_\mx - \vec g\right)
+ \nabla \beta_\mx(S^\ve_\mx) \bigg\}.
$$
By using the estimate \eqref{beta-un} and the assumptions {(A.2)}, {(A.4)},
we get the uniform bound:
\begin{equation}
\label{Fve}
\left\Vert {\EuScript F}^{\,\ve} \right\Vert_{L^2(\Omega^\ve_{\mx,T})} \leqslant C.
\end{equation}

Let define a function:
\begin{equation*}
\varphi_w\left(x, \frac{x}{\ve}, t \right) \in {\EuScript D}(\Omega_T; C^\infty_{per}(Y))
\quad {\rm such\,\, that\,\,} \varphi_w = 0 \,\,\, {\rm for}\,\,\, y\in Y_\fr.
\end{equation*}
Plugging $\varphi_w$ in (\ref{wf-1-gl-GP}), and taking into account condition {(A.9)},
we get:
\begin{equation}
\label{eqn}
-\int\limits_{\Omega_{T}} {\bf 1}^\ve_\mx(x)\, \Phi^\ve_\mx(x)\, S^\ve_\mx\,
\frac{\partial \varphi_w}{\partial t} \, dx\, dt
+
\ve^{\frac{\theta}{2}}\, \int\limits_{\Omega^\ve_{\mx,T}} {\EuScript F}^{\,\ve} \, \nabla_x \varphi_w \, dx\, dt
+
\ve^{\frac{\theta}{2}-1}\, \int\limits_{\Omega^\ve_{\mx,T}} {\EuScript F}^{\,\ve} \, \nabla_y \varphi_w \, dx\, dt = 0,
\end{equation}
We pass to the two-scale limit in \eqref{eqn} using \eqref{Fve}. We obtain:
\begin{equation}
\label{eqn-1}
\int\limits_{\Omega_{T}\times Y_\mx} \Phi_\mx (y)\, s(x, y, t)\,
\frac{\partial \varphi_w}{\partial t} \, dx\, dt\, dy = 0.
\end{equation}
This completes the proof of Lemma \ref{nosource}. \fin

Finally, from the equations (\ref{1newpoint-7:2ok}), (\ref{hernya-18ok}) in view of
Lemma \ref{nosource}, arguing as in subsection \ref{hom-eq-in-ph-pres}, we arrive to
the desired system (\ref{H-0-veryhigh}). This completes the proof of
Theorem \ref{t-hom-main>2}. \fin

\subsection{Moderate contrast media: ${\bf 0<}{\boldsymbol \theta}{\bf<2}$}
\label{moder-case-subsec}

We study the asymptotic behavior of the solution to problem
(\ref{debut2}) as $\ve \to 0$ in the case $\varkappa(\ve) = \ve^\theta$
with $0 < \theta < 2$. In particular, we are going to show that the effective model reads:
\begin{equation}
\label{H-0-moderate}
\left\{
\begin{array}[c]{ll}
0 \leqslant S \leqslant 1 \quad {\rm in} \,\, \Omega_T; \\[2mm]
\displaystyle
\frac{\partial }{\partial t}
\left[\Phi^\star(x)\, S + \widehat \Phi_\mx\, {\EuScript P}(S) \right]
- {\rm div}_x\, \bigg\{\mathbb{K}^\star(x)\, \lambda_{\,\fr,w}(S) \big(\nabla P_w - \vec g \big)  \bigg\}
= F^\star_w \quad {\rm in} \,\, \Omega_T;
\\[5mm]
\displaystyle
- \frac{\partial }{\partial t}
\left[\Phi^\star(x)\, S + \widehat \Phi_\mx\, {\EuScript P}(S) \right]
- {\rm div}_x\, \bigg\{\mathbb{K}^\star(x)\, \lambda_{\,\fr,n}(S)
\big(\nabla P_n - \vec g \big)  \bigg\}
= F^\star_n \quad {\rm in} \,\, \Omega_T;\\[5mm]
P_{\fr,c}(S) = P_n - P_w \quad {\rm in} \,\, \Omega_T,
\end{array}
\right.
\end{equation}
where the effective porosity $\Phi^\star$, the effective source terms $F^\star_w, F^\star_n$,
and the homogenized permeability tensor $\mathbb{K}^\star$
in \eqref{H-0-moderate} are defined in \eqref{H-1}, \eqref{H-1FFF} and
\eqref{H-2}, respectively.

The boundary conditions and the initial conditions
for the system (\ref{H-0-moderate}) are given by (\ref{H-7}), (\ref{H-9}).

In this case we observe a complete decoupling between microscale and macroscale, which
is not the case for the critical scaling $\theta = 2$.

The third main result of the paper is as follows.

\begin{theorem}
\label{t-hom-main<2}
Let $\varkappa(\ve) = \ve^\theta$ with $0 < \theta < 2$ and
let assumptions {(A.1)-(A.9)} be fulfilled. Then the solution of the {\rm initial problem
(\ref{debut2}), (\ref{inter-condit})-(\ref{init1})} converges (up to a subsequence) in the two-scale
sense to a weak solution of the {\rm homogenized problem \eqref{H-0-moderate},
(\ref{H-7}), (\ref{H-9})}.
\end{theorem}

Let $0 < \theta < 2$. In the proof of Theorem \ref{t-hom-main<2} we follow
the lines of the proof of Theorem \ref{t-hom-main}. Namely, arguing as
in Sections \ref{passage1}, \ref{passage2}, \ref{ident-p-beta}, we obtain
the homogenized equations (\ref{1newpoint-7:2ok}), (\ref{hernya-18ok}).
Namely, in the case of the moderate contrast we have:
\begin{equation}
\label{H-0-moderate-proof}
\left\{
\begin{array}[c]{ll}
0 \leqslant S \leqslant 1 \quad {\rm in} \,\, \Omega_T; \\[2mm]
\displaystyle
\Phi^\star(x)\, \frac{\partial S}{\partial t}
- {\rm div}_x\, \bigg\{\mathbb{K}^\star(x)\, \lambda_{\,\fr,w}(S) \big(\nabla P_w - \vec g \big)  \bigg\}
=
\widehat{{\EuScript Q}}_w + F^\star_w \quad {\rm in} \,\, \Omega_T;
\\[5mm]
\displaystyle
- \Phi^\star(x)\, \frac{\partial S}{\partial t}
- {\rm div}_x\, \bigg\{\mathbb{K}^\star(x)\, \lambda_{\,\fr,n}(S)
\big(\nabla P_n - \vec g \big)  \bigg\}
=
\widehat{{\EuScript Q}}_n + F^\star_n \quad {\rm in} \,\, \Omega_T;\\[5mm]
P_{\fr,c}(S) = P_n - P_w \quad {\rm in} \,\, \Omega_T,
\end{array}
\right.
\end{equation}
where the effective porosity $\Phi^\star$, the effective source terms $F^\star_w, F^\star_n$,
and the homogenized permeability tensor $\mathbb{K}^\star$
in \eqref{H-0-moderate-proof} are defined in \eqref{H-1}, \eqref{H-1FFF} and
\eqref{H-2}, respectively. For any $x \in \Omega$ and $t > 0$, the matrix-fracture source terms
$\widehat{{\EuScript Q}}_w$, $\widehat{{\EuScript Q}}_n$ in \eqref{H-0-moderate} have the form:
$$
\widehat{{\EuScript Q}}_w \eqdef - \widehat \Phi_\mx\,\frac{\partial s}{\partial t}(x, t)
= - \widehat{{\EuScript Q}}_n \quad {\rm with} \,\,\,
\widehat \Phi_\mx \eqdef \frac{1}{|Y_\mx|}\, \int\limits_{Y_\mx} \Phi_\mx(y)\,dy.
$$

In order to complete the proof of Theorem \ref{t-hom-main<2}, we have to
identify  the saturation function $s$ appearing on the right-hand side of
equations in (\ref{H-0-moderate-proof}). The following result holds true:

\begin{lemma}
\label{our-m3as-version}
Let $s$ be the weak limit of $\{\mathfrak{D}^\ve S_{\mx}^\ve\}_{\ve>0}$ and
$S$ is the saturation function defined in (\ref{2s-1}). Then
\begin{equation}
\label{+0-3}
s = {\EuScript P}(S) \quad {\rm a.e.\,\, in}\,\, \Omega_T
\quad {\rm with} \,\, {\EuScript P}(S) = (P_{c,m}^{-1} \circ P_{c,f})(S).
\end{equation}
\end{lemma}

\noindent{\bf Proof of Lemma \ref{our-m3as-version}}. Applying Lemma~\ref{lem-hitro-vyeb}
and Proposition~\ref{prop-vyeb} we conclude that, for any
$x_0 \in \Omega\setminus {\EuScript A}_{\bf n}$,
\begin{align*}
\beta_\mx(s^\ve_{\mx,x_0}) &\to \beta_\mx(s_{x_0})\quad \text{weakly in }\;
L^2(0,T; H^1(Y_{\mx})),\\
s^\ve_{\mx,x_0} &\to s_{x_0} \quad \text{a.e. in }\; Y_{\mx}\times (0,T),
\end{align*}
and the limit $s_{x_0}$ does not depend of the fast variable $y$. Due to
continuity of the trace operator we also have:
\begin{align*}
\beta_\mx(s^\ve_{\mx,x_0})\big|_{\Gamma_{\mx\fr}} \to
\beta_\mx(s_{x_0})\big|_{\Gamma_{\mx\fr}}
\quad \text{weakly in }\; L^2(0,T; L^2(\Gamma_{\mx\fr})).
\end{align*}
On the other hand we know that, for a.e. $x_0\in \Omega$,
\begin{align*}
{\cal M}(\beta_\fr(\mathfrak{D}^\ve(\widetilde S^\ve_\fr(x_0,\cdot,\cdot))))
\big|_{\Gamma_{\mx\fr}} = \beta_{\mx}(s^\ve_{\mx,x_0})
\big|_{\Gamma_{\mx\fr}} \quad {\rm with} \,\,\,
{\cal M} \eqdef \beta_{\mx} \circ ( P_{\mx,c})^{-1}\circ P_{\fr,c} \circ (\beta_{\fr})^{-1}
\end{align*}
a.e. on $\Gamma_{\mx\fr}\times (0,T)$.
For a.e. $x_0\in \Omega$, from Corollary~\ref{corol-dilop} we have that
\begin{align*}
{\cal M}(\beta_\fr(\mathfrak{D}^\ve(\widetilde S^\ve_\fr(x_0,\cdot,\cdot)))) \to
{\cal M}(\beta_\fr(S(x_0,\cdot,\cdot)))\quad\text{strongly in }\;
L^2(0,T; L^2(\Gamma_{\mx\fr}))
\end{align*}
and therefore, for a.e. $x_0 \in \Omega\setminus {\EuScript A}_{\bf n}$,
\begin{align*}
\beta_\mx(s_{x_0})\big|_{\Gamma_{\mx\fr}} =
{\cal M}(\beta_\fr(S(x_0,\cdot,\cdot)))\big|_{\Gamma_{\mx\fr}}.
\end{align*}
Since these functions are independent of $y$ we have that $\beta_\mx(s_{x_0}) =
{\cal M}(\beta_\fr(S(x_0,\cdot))$ in $L^2(0,T)$, or, equivalently,
$s_{x_0} = {\EuScript P}(S(x_0,\cdot))$. Now, for a chosen
$x_0 \in \Omega\setminus {\EuScript A}_{\bf n}$, we can find a subsequence
such that
\begin{align*}
s^\ve_{\mx,x_0} &\to {\EuScript P}(S(x_0,\cdot)) \quad \text{a.e. in }\;
Y_{\mx}\times (0,T).
\end{align*}
Since the limit is uniquely defined  by the limit $S$ of the sequence $\mathfrak{D}^\ve(\widetilde S^\ve_\fr)$
we conclude that the whole sequence converge to the same limit (that is the whole subsequence for which $\mathfrak{D}^\ve(\widetilde S^\ve_\fr)$
converges). Now we can repeat our procedure for almost
any  $x_0 \in \Omega\setminus {\EuScript A}_{\bf n}$ and conclude that
$s = {\EuScript P}(S)$ a.e. in $(\Omega\setminus {\EuScript A}_{\bf n})\times (0,T)$.
Thanks to Propositions  \ref{proppi-1}, \ref{proppi-2}, the measure of the set ${\EuScript A}_{\bf n}$
goes to zero as ${\bf n}\to\infty$ and the desired equality (\ref{+0-3}) is proved. \fin

Now we complete easily the proof of Theorem \ref{t-hom-main<2}. Taking into account
(\ref{+0-3}) we can rewrite (\ref{H-0-moderate-proof}) and thereby obtain (\ref{H-0-moderate}).
Theorem \ref{t-hom-main<2} is proved. \fin

\section*{Acknowledgments}
The work of M. Jurak and A. Vrba\v ski was funded by
Croatian science foundation project no 3955.  The work of L. Pankratov was funded by the
RScF research project N 15-11-00015.
This work was partially supported by ISIFoR (http://www.carnot-isifor.eu/) France. The
supports are gratefully acknowledged. Most of the work on this paper was done when M. Jurak and L. Pankratov were visiting the Applied Mathematics Laboratory of the University of Pau \& CNRS.

\renewcommand{\baselinestretch}{0.6}

\begin{small}

\end{small}

\end{document}